\documentclass[a4paper,psamsfonts]{amsart}
\usepackage{amssymb,mathrsfs,pstricks,pst-node}
\newtheorem{lem}{Lemma}
\newtheorem{thm}[lem]{Theorem}
\newtheorem{defn}[lem]{Definition}
\newtheorem{cor}[lem]{Corollary}
\newtheorem{prop}[lem]{Proposition}

\DeclareMathOperator{\Stab}{Stab}
\DeclareMathOperator{\grp}{gp}

\DeclareMathOperator{\Ref}{Ref}

\DeclareMathOperator{\Odd}{Odd}
\DeclareMathOperator{\EOdd}{EOdd}

\DeclareMathOperator{\FC}{FC}
\DeclareMathOperator{\Even}{E}
\DeclareMathOperator{\oddgraph}{\Omega}
\newcommand\Z{\mathbb{Z}}
\newcommand\R{\mathbb{R}}
\numberwithin{equation}{section}
\let\le\leqslant
\let\ge\geqslant
\begin{document}

\title{Reflections in abstract Coxeter groups}

\author{W.~N.~Franzsen}
\address{Australian Catholic University\\
25A Barker Rd\\
Strathfield NSW 2135\\
Australia}
\email{b.franzsen@mary.acu.edu.au}

\author{R.~B.~Howlett}
\address{School of Mathematics and Statistics\\
University of Sydney\\
NSW 2006\\ Australia}
\email{R.Howlett@maths.usyd.edu.au}

\author{B.~M\"uhlherr}
\address{D\'epartement de Math\'ematiques\\
ULB C.P. 216\\ Bd. du Triomphe \\ 1050 Bruxelles \\
Belgium}
\email{bernhard.muhlher@ulb.ac.be}

\abstract
Let $W$ be a Coxeter group and $r\in W$ a reflection. If the group of
order 2 generated by $r$ is the intersection of all the maximal
finite subgroups of $W$ that contain it, then any isomorphism from
$W$ to a Coxeter group $W'$ must take $r$ to a reflection in~$W'$.
The aim of this paper is to show how to determine, by inspection of
the Coxeter graph, the intersection of the maximal finite sugroups
containing~$r$. In particular we show that the condition above is satisfied
whenever $W$ is infinite and irreducible, and has the property that all
rank two parabolic subgroups are finite. So in this case all
isomorphisms map reflections to reflections.
\endabstract

\maketitle
\section{Introduction}

The dihedral group of order 12 can be considered as Coxeter
group of type $I_2(6)$ or as Coxeter group  of type $A_1 \times I_2(3)$.
This example shows that, in general, the set of reflections in
a Coxeter system is not determined by the abstract group $W$
alone, but does indeed depend on the choice of the Coxeter generating
set $R$. However there are a lot of examples of Coxeter systems
$(W,R)$ where the abstract group does determine the set of
reflections or even the set $R$ up to $W$-conjugacy. The main
motivation for the present paper is to show that 
the latter holds for infinite
Coxeter groups having a finite, irreducible and 2-spherical
Coxeter generating set, which is our Theorem \ref{thm2sph} below.

In view of the main result of \cite{CM} it suffices to show that
these Coxeter groups determine the set of reflections. In order to
achieve this goal we provide a handy criterion for an involution
in an abstract Coxeter group $W$ to be a reflection with respect to 
any Coxeter generating set of $W$. Our principal observation is the
following. Let $(W,R)$ be a Coxeter system and let $w \in W$
be an involution. If $w \not \in R^W$, then the centralizer
of $w$ in $W$ contains  a finite normal subgroup properly containing
$\langle w \rangle$. This is an immediate consequence of
Richardson's result in \cite{RR}. Thus, if $w \in W$
is an involution having the property that $\langle w \rangle$
is a maximal finite normal subgroup of its centralizer in $W$,
then $w$ is a reflection with respect to any Coxeter generating
set of $W$.

It turns out that it is more convenient to work with
the \textit{finite continuation\/} of an involution rather than
to consider finite normal subgroups of its centralizer. 
The finite continuation of a finite
order element $w$ in a Coxeter group is defined to be the intersection
of all maximal finite subgroups containing it; we write
$\FC(w)$ for the finite continuation of~$w$. In this paper we
restrict our attention to finitely generated Coxeter groups.
For these it is a consequence of a result of Tits 
that every element of finite order is contained in some maximal
finite subgroup; so $\FC(w)$ is a finite subgroup of $W$
(see Corollary \ref{FCexis}
below). The main
result of the present paper is a complete description of 
the finite continuation of a simple reflection in a Coxeter system
of finite rank. Its proof constitutes the bulk of this paper.

\smallskip
\noindent
{\bf Main Result:} {\it Let $(W,R)$ be a Coxeter system of finite rank. 
Then the following holds.

\begin{itemize}
\item[a)] For each $r \in R$ the finite continuation of $r$ can
be described.
\item[b)] Given an involution $w \in W$ such that 
$\FC(w) = \langle w \rangle$, then $w \in R^W$.
\end{itemize}
}

\smallskip
Part a) of our main result is Theorem \ref{main}. Its precise statement
requires some preparation.  Part b) is Corollary \ref{FCcrit1}.

The main result of this paper is in fact the first of two steps
to reduce the isomorphism problem for Coxeter groups
to its `reflection-preserving' version. The second step is
given in \cite{HM}. We refer to \cite{toro} for further information
about the applications to the general isomorphism problem.
 
A special instance of the isomorphism problem for Coxeter groups
is the question about their rigidity (see \cite{BMMN} for further
information). In combination with the main result of 
\cite{CM} a consequence of
 our main result is the following rigidity result.

\begin{thm} \label{thm2sph}
Let $(W,R)$ be an irreducible, non-spherical Coxeter system such that
$R$ is finite and such that $rr'$ has finite order for all $r,r' \in R$.
Then the following assertions hold.

\begin{itemize}
\item[a)] For each $r \in R$ we have $\FC(r) = \langle r \rangle$.
\item[b)] If $S \subseteq W$ is such that $(W,S)$ is a Coxeter system,
then there exists $w \in W$ such that $S^w = R$.
\item[c)] All automorphisms of $W$ are inner-by-graph.
\end{itemize}

\end{thm}
In the language of \cite{BMMN}, Part b) of the previous theorem means
that an infinite, irreducible, 2-spherical Coxeter system is 
strongly rigid. Part c), which is an immediate consequence of Part b),
improves the result of \cite{HRT}.

To conclude this introduction we remark that characterizations results
for reflections in even Coxeter groups have been obtained in \cite{BM}.
Some of the results there can be deduced as corollaries of our 
main result as well.

\subsection*{Acknowledgement}
The authors thank Fr\'ed\'eric Haglund for a helpful discussion
on the subject.

\section{Precise Statement of the Main Result}

Recall that a Coxeter group is a group with a presentation of the
form
\begin{equation}\label{eq:pres}
W=\grp\langle\,\{\,r_a\mid a\in\Pi\,\}\mid(r_ar_b)^{m_{ab}}=1
\text{ for all }a,b\in\Pi\,\rangle
\end{equation}
where $\Pi$ is some indexing set, whose cardinality is called the
\textit{rank\/} of $W$ (relative to this presentation), and the
$m_{ab}$ satisfy the following conditions: $m_{ab}=m_{ba}$, each
$m_{ab}$ lies in the set
$\{\,m\in\mathbb{Z}\mid m\ge 1\,\}\cup\{\infty\}$, and $m_{ab}=1$
if and only if $a=b$. When $m_{ab}=\infty$ the relation
$(r_ar_b)^{m_{ab}}=1$ is interpreted as vacuous.
We shall restrict attention to Coxeter groups of finite rank.

A \textit{reduced expression\/} for an element $w\in W$ is a minimal
length word expressing $w$ as a product of elements of the
distinguished generating set $\{\,r_a\mid a\in\Pi\,\}$. We define
$\ell(w)$ to be the length of a reduced expression for~$w$.

As is well known (and as we shall describe in Section~3 below),
every Coxeter group W can be realized geometrically as a group
generated by reflections. In this realization of~$W$ the reflections in
$W$ are the conjugates of the generators~$r_a$.

The \textit{Coxeter graph\/} associated with the presentation above is
the graph with vertex set $\Pi$ and edge set consisting of those pairs of
vertices $\{a,b\}$ for which $m_{ab}\ge 3$. The edge $\{a,b\}$ is given the
label~$m_{ab}$. The \textit{components\/} of $\Pi$ are the connected
components of the graph, and we say that $W$ is \textit{irreducible\/}
if the graph is connected.

For each $I\subseteq\Pi$ we define $W_I$ to be the subgroup of $W$
generated by the set $\{\,r_a\mid a\in I\,\}$; we call these subgroups
the \textit{visible\/} subgroups of~$W$. A \textit{parabolic subgroup\/}
of $W$ is any conjugate of a visible subgroup. We say that $I\subseteq\Pi$
is \textit{spherical\/} if $W_I$ is finite, and we say that $\Pi$ (or
$W$) is \textit{$k$-spherical} if all $k$-element subsets of $\Pi$
are spherical.

The defitions given so far are fairly standard. In order to
facilitate
the precise statement
of the main result, we introduce some nonstandard notation and terminology
(in Definitions~\ref{FC}, \ref{EO}, \ref{c3neighbour}, \ref{focus}
and \ref{halffocus} below).

\begin{defn}\label{FC} {\rm If $w\in W$ has finite order, define the
\textit{finite continuation\/} of $w$, written $\FC(w)$, to be the
intersection of all maximal finite subgroups of~$W$ containing~$w$.}
\end{defn}

\begin{defn}\label{EO} {\rm The \textit{odd graph\/} of $W$ is the graph
$\oddgraph(\Pi)$ obtained from the Coxeter graph by deleting the edges whose
labels are infinite or even. For each $a\in\Pi$ we  define $\Odd(a)$ to
be the
connected component of $\oddgraph(\Pi)$ containing~$a$. For each connected
component $M$ of $\oddgraph(\Pi)$ we define $\Even(M)$ to be the union
of $M$ with the set of all $b\in\Pi$ such that $m_{cb}$ is even for some
$c\in M$. We also abbreviate $\Even(\Odd(a))$ to $\EOdd(a)$.}
\end{defn}

In the discussions below, when we refer to the components of $\Even(M)$
we regard $\Even(M)$ as the full subgraph of the Coxeter graph
spanned by the vertices in~$\Even(M)$. In other words, the edges with even
and infinite labels, deleted when forming the odd graph, are restored
in $\Even(M)$.

Note that if $a\in L\subseteq\Pi$ and $W_L$ is finite then
$m_{ab}<\infty$ for all $b\in L$. Whether $m_{ab}$ is odd or even
it follows that $b\in\EOdd(a)$. Thus $L\subseteq\EOdd(a)$.

\begin{defn}\label{c3neighbour}
{\rm Let $M\subseteq\Pi$ be a connected component of $\oddgraph(\Pi)$.
We call $b\in\Pi\setminus M$ a $C_3$-\textit{neighbour\/} of $M$ if
$m_{bc}\in\{2,4\}$ for all $c\in E(M)$, the case $m_{bc}=4$ occurring
for at least one~$c$, and for each $c\in E(M)$ with $m_{bc}=4$ there
is an $a\in M$ such that the following conditions are satisfied:
\begin{enumerate}
\item $m_{ba}=2$ and $m_{ca}=3$, and $m_{cd}=\infty$ for all
$d\in M\setminus\{a,c\}$;
\item for all $e\in\Pi\setminus(M\cup\{b\})$, either $m_{ce}=\infty$ or
$m_{ae}=m_{ce}=m_{be}=2$.
\end{enumerate}}
\end{defn}

\begin{defn}\label{focus}
{\rm Let $M\subseteq\Pi$ be a connected component of $\oddgraph(\Pi)$, and let
$a\in M$ and $b\in\Pi\setminus M$.
We call the pair $(a,b)$ a \textit{focus\/} of $M$  in $\Pi$ if the
following conditions all hold.
\begin{enumerate}
\item
All the edge labels of $M$ are~3, and $M$ is a tree.
\item
For each $c\in M$, the set $C[b..c]\subseteq\Pi$ consisting of $b$ and those
elements of $M$ that form the path from $a$~to~$c$ in~$M$ constitutes a system
of type~$C_k$ (for some~$k$ dependent on~$c$).
\item
If $c,\,d\in M\cup\{b\}$ with $c\notin C[b..d]$ and $d\notin C[b..c]$ then
$m_{cd}=\infty$.
\item
If $m_{ce}\ne\infty$ for some $c\in M$ and $e\in\Pi\setminus(M\cup\{b\})$,
then $m_{ce}=2=m_{de}$ for all $d\in C[b..c]$.
\item
The vertices of $M\cup\{b\}$ do not form a spherical component of $\Even(M)$.
\end{enumerate}}
\end{defn}

\begin{defn}\label{halffocus}
{\rm Let $M\subseteq\Pi$ be a connected component of $\oddgraph(\Pi)$, and let
$a,\,b\in M$. We call the two-element set $\{a,b\}$ a \textit{half focus\/}
of $M$ in $\Pi$ if $m_{ab}=2$ and the following conditions all hold.
\begin{enumerate}
\item
We have $m_{ac}=m_{bc}\in\{2,3\}$ for all $c\in M\setminus\{a,b\}$, and
$m_{ac}=m_{bc}\in\{2,\infty\}$ for all $c\in\Pi\setminus M$.
\item
All the edge labels of $M\setminus\{b\}$ are~3, and $M\setminus\{b\}$
is a tree.
\item
For each $c\in M\setminus\{a,b\}$, the set $D[a,b..c]\subseteq\Pi$
consisting of $b$ and those and those elements of $M\setminus\{b\}$ that
form the path from $a$~to~$c$ constitutes a system
of type~$D_k$ (for some~$k$ dependent on~$c$).
\item
If $c,\,d\in M\setminus\{a,b\}$ with $c\notin D[a,b..d]$ and $d\notin D[a,b..c]$
then $m_{cd}=\infty$.
\item
If $m_{ce}\ne\infty$ for some $c\in M\setminus\{a,b\}$ and
$e\in\Pi\setminus M$, then $m_{ce}=2=m_{de}$ for all $d\in D[a,b..c]$.
\item
The vertices of $M$ do not form a spherical component of $\Even(M)$.
\end{enumerate}}
\end{defn}

We are now able to give a precise statement of Part a) of our main result.

\begin{thm}\label{main}
For each connected component $M$ of $\oddgraph(\Pi)$ there is at least
one $a\in M$ such that $\FC(r_a)$ is a visible subgroup of~$W$. We have
the following possibilities.

\smallskip
\noindent{\rm Case A:}
Suppose that the component of $\Even(M)$ containing~$M$ is spherical,
and let $a\in M$ be arbitrary. Then $\FC(r_a)=W_J$, where $J$ is the
union of the spherical components of~$\Even(M)$.

\smallskip
\noindent{\rm Case B:}
Suppose that the component of $\Even(M)$ containing~$M$ is not spherical,
and $M$ does not have any focus or half-focus in~$\Pi$, and let $J'$ be the
union of the spherical components of $\Even(M)$ and the set of $C_3$-neighbours
of~$M$. If $a\in M$ is not adjacent in $\Pi$ to any $C_3$-neighbour of~$M$
then $\FC(r_a)=W_{J'\cup\{a\}}$, and if $a\in M$ is adjacent in $\Pi$ to a
$C_3$-neighbour of~$M$ then $\FC(r_a)$ is not visible.

\smallskip
\noindent{\rm Case C:}
Suppose that $(a,b)$ is a focus of~$M$. Then $\FC(r_a)=W_J$
where $J$ is the union of $\{a,b\}$ and the spherical components of $\Even(M)$,
and $\FC(r_c)$ is not visible for any $c\in M\setminus\{a\}$.

\smallskip
\noindent{\rm Case D:}
Suppose that $\{a,b\}$ is a half-focus of~$M$. Then $\FC(r_a)=\FC(r_b)=W_J$,
where $J$ is the union of $\{a,b\}$ and the spherical components of $\Even(M)$,
and $\FC(r_c)$ is not visible for any $c\in M\setminus\{a,b\}$.
\end{thm}

\section{Reflections and root systems}

Let $\R$ be the real field, and $V$ the vector space over $\R$ with
basis~$\Pi$. Let $B$ the bilinear form on $V$ such that for all
$a,b\in\Pi$,
\[
B(a,b)=-\cos(\pi/m_{ab}).
\]
To make our notation more compact we define $u\cdot v=B(u,v)$ for all
$u,v\in V$. Note that $a\cdot a=1$ for all $a\in\Pi$, since
$m_{aa}=1$.

For each $a\in V$ such that $a\cdot a=1$, the \textit{reflection
along~$a$} is the transformation of $V$ given by $v\mapsto v-2(a\cdot v)a$.
It is well known (see, for example, Corollary 5.4 of \cite{JH}) that $W$
has a faithful representation on $V$ such that, for all $a\in\Pi$, the
element $r_a$ acts as the reflection along~$a$. We shall identify elements
of $W$ with their images in this
representation. We also use the notation $r_a$ for the reflection
along~$a$ whenever $a\in V$ satisfies $a\cdot a=1$. It is
straightforward to show that each reflection~$r_a$ preserves the
form~$B$; hence all elements of $W$ preserve~$B$. Furthermore, the
equation $gr_ag^{-1}=r_{ga}$ holds for all $a\in V$ such that
$a\cdot a=1$ and all transformations $g$ that preserve~$B$.

We write $\Ref(W)$ for the set of all reflections in $W$. It is
immediate from the above comments that if
$\Phi=\{\,w a\mid w\in W,\ a\in \Pi\,\}$ then
$\{\,r_b\mid b\in\Phi\,\}\subseteq\Ref(W)$.

The set $\Phi$ is called the \textit{root system\/} of $W$, and
elements of $\Phi$ are called \textit{roots}. Elements of the basis
$\Pi$ are called \textit{simple roots}, and the reflections $r_a$ for
$a\in \Pi$ are called \textit{simple reflections}. A root is said to
be \textit{positive\/} if it has the form $\sum_{a\in\Pi}\lambda_a a$
with $\lambda_a\ge 0$ for all $a\in\Pi$, and \textit{negative\/}
otherwise. We write $\Phi^+$ for the set of all positive roots and
$\Phi^-$ for the set of all negative roots.

\begin{lem} \label{somebasicfacts}
With the notation as above, the following statements hold.
\begin{enumerate}
\item Every negative root has the form $\sum_{a\in \Pi}\lambda_a a$
with $\lambda_a\le 0$ for all $a\in\Pi$. Furthermore,
$\Phi^-=\{\,-b\mid b\in\Phi^+\,\}$.
\item If $w\in W$ and $a\in\Pi$ then
\[
\ell(wr_a)=
\begin{cases}
\ell(w)+1&\text{if $wa\in\Phi^+$,}\\
\ell(w)-1&\text{if $wa\in\Phi^-$.}
\end{cases}
\]
\item If $t\in\Ref(W)$ then $t=r_b$ for some $b\in \Phi$.
\item The group $W$ is finite if and only if the bilinear form $B$ is
positive definite.
\item The root system $\Phi$ is finite if and only if the group $W$
is finite.
\end{enumerate}
\end{lem}

\begin{proof} Proofs of (1) and (2) can be found in
\cite[Section 5.4]{JH}, Theorem 4.1 in \cite{VD} includes both
(4)~and~(5), and (3) is \cite[Lemma 2.2]{HRT}.
\end{proof}

The following result is well known.

\begin{lem}\label{class} Let $a\in\Pi$. Then $\Odd(a)=\Pi\cap Wa$.
\end{lem}

For each $w\in W$ we define $N(w)=\{\,b\in\Phi^+\mid wb\in\Phi^-\,\}$.
By Part (2) of Lemma~\ref{somebasicfacts}, if $w\ne 1$ then
$N(w)\cap\Pi\ne\emptyset$. An easy induction shows that $N(w)$ has
exactly $\ell(w)$ elements. In particular, $N(w)$ is a finite set. It
is also easily shown that if $\Phi$ is finite then there is a unique
$w\in W$ such that $N(w)=\Phi^+$. This element, which we denote by
$w_\Pi$, is also the unique element of maximal length in~$W$ (which is
a finite group). Furthermore, $w_\Pi\Pi=-\Pi$.

For each $\Gamma\subseteq\Phi$ the subgroup $W_\Gamma$ generated by the
set $\{\,r_a\mid a\in\Gamma\,\}$ is called a \textit{reflection subgroup\/}
of~$W$. The set $\Phi_\Gamma=\{\,a\in\Phi\mid r_a\in W_\Gamma\,\}$ is called
the \textit{root subsystem\/} generated by~$\Gamma$. Let
$\Phi_\Gamma^+=\Phi_\Gamma\cap\Phi^+$ and $\Phi_\Gamma^-=\Phi_\Gamma\cap\Phi^-$,
and define
$$
\Pi_\Gamma=\{\,a\in\Phi_\Gamma^+\mid N(r_a)\cap\Phi_\Gamma=\{a\}\,\}.
$$
The main theorem of Deodhar \cite{VD2} and Theorem~(3.3) of Dyer~\cite{MD}
yield the following result.

\begin{thm}\label{base}
For each $\Gamma\subseteq\Phi$ the group $W_\Gamma$ is a Coxeter
group on the generating set $\{\,r_a\mid a\in\Pi_\Gamma\,\}$. The set
$\{\,a\cdot b\mid a,\,b\in\Pi_\Gamma\text{ and }a\ne b\,\}$ is a subset
of $\mathscr C=\{\,-\cos(\pi/m)\mid 2\le m\in\Z\,\}\cup(-\infty,-1]$. Moreover,
if $\Delta$ is any subset of $\Phi^+$ such that
$\{\,a\cdot b\mid a,\,b\in\Delta\text{ and }a\ne b\,\}\subseteq \mathscr{C}$
then $W_\Delta$ is a Coxeter group on the generating set
$\{\,r_a\mid a\in\Delta\,\}$.
\end{thm}

Note that the notation $W_\Gamma$ introduced above is an extension of
the notation for visible subgroups introduced in Section~2. However,
if $\Gamma\nsubseteq\Pi$ then $W_\Gamma$ need not be visible.

It is clear that if $I\subseteq\Pi$ then $W_I$ preserves the subspace
$V_I$ of $V$ spanned by $I$, and acts on this subspace as a Coxeter
group with $I$ as its set of simple roots. In this case $\Phi_I=\Phi\cap V_I$
and $\Pi_I=I$.

The following simple facts are well known.

\begin{lem}\label{subrootsys}
In the above situation, $\Phi_I=\Phi\cap V_I$. Furthermore,
$w\in W$ normalizes $W_I$ if and only if $w\Phi_I=\Phi_I$.
In particular, for all $a\in\Phi$, the reflection $r_a$
normalizes $W_I$ if and only if $a\in\Phi_I$ or $a\cdot b=0$ for
all~$b\in I$.
\end{lem}

Suppose that $I\subseteq\Pi$ and $a\in\Pi\setminus I$, and let
$L$ be the component of (the Coxeter graph of)
$I\cup\{a\}$ to which $a$ belongs. If $W_L$ is finite we define
$v[a,I]=w_Lw_{L\setminus\{a\}}$. It is easily seen that
$v[a,I]I\subseteq I\cup\{a\}$, and that $v[a,I]b=b$ for all
$b\in I\setminus L$. In particular,
$v[a,I]I\in\mathscr{I}=\{\,J\subseteq\Pi\mid J=wI
\text{ for some $w\in W$ }\}$.
It was proved in \cite{H} (for finite Coxeter
groups) and in \cite{VD} (in the general case) that every element
$w\in W$ satisfying $wI\subseteq\Pi$ can be expressed as a product
of elements of the form $v[a,I']$, with $I'\in\mathscr{I}$ and
$a\in\Pi\setminus I'$. That is,
\begin{equation}\label{canonexp}
w=v[a_1,I_1]v[a_2,I_2]\cdots v[a_n,I_n]
\end{equation}
for some $I_i,\,a_i$ such that (for each $i$) the component
of $I_i\cup\{a_i\}$ containing~$a_i$ corresponds to a finite visible
subgroup, $v[a_i,I_i]I_i=I_{i-1}$ for $1<i\le n$, and $I_n=I$.
Furthermore, the following result holds.

\begin{prop}\label{normalizers} Let $I,\,J\subseteq\Pi$. Then
$\{\,w\in W\mid wW_Iw^{-1}=W_J\,\}=N(J,I)W_I$, where
$N(J,I)=\{\,w\in W\mid wI=J\,\}$. Furthermore, for each $w\in N(J,I)$
and each $a\in\Pi\cap N(w)$ there is an expression for $w$ of the form
(\ref{canonexp}) above, with $(a_n,I_n)=(a,I)$ and
$\ell(w)=\sum_{i=1}^n\ell(v[a_i,I_i])$.
\end{prop}

The following lemma, which appears in
\cite[Exercise 2d, p.~130]{NB}, is fundamental to all of our
arguments.

\begin{lem}[Tits] \label{lemtits}
If\/ $W$ is a Coxeter group and $H\le W$ is finite, then $H$ is
contained in a finite parabolic subgroup of $W$.
\end{lem}

One immediate consequence of Lemma~\ref{lemtits} is that every maximal
finite subgroup of a Coxeter group is parabolic. Another consequence
of the previous lemma 
is
that each finite subgroup of $W$ is contained in a maximal finite parabolic
subgroup. (Remember that we always assume that $W$ is finitely generated.) 
Thus the set of maximal finite subgroups of $W$ containing
a given finite order element of $W$ is not empty, and hence
we have the following fact.

\begin{cor} \label{FCexis}
If $w \in W$ has finite order, then $\FC(w)$ is a well-defined finite
subgroup of $W$.
\end{cor}

\begin{lem}[Kilmoyer] \label{kilmoyerlem}
Let $I,J\subseteq \Pi$. Then every $(W_I,W_J)$ double coset in $W$
contains a unique element of minimal length; moreover, if $d$ is the
minimal length element of\/ $W_I d W_J$ then $W_I\cap dW_Jd^{-1}=W_K$,
where $K=I\cap d J$.
\end{lem}

\begin{proof}
See \cite[Theorem 2.7.4]{RC}.
\end{proof}

\begin{cor} \label{intersectionofparasispara}
The intersection of a finite number of parabolic subgroups is a
parabolic subgroup.
\end{cor}

The following consequence of Lemmas \ref{lemtits}~and~\ref{kilmoyerlem}
is proved in \cite[Lemma 11]{FH2}.

\begin{lem} \label{maxvismax}If $J$ is a maximal spherical subset of $\Pi$
then $W_J$ is a maximal finite subgroup of~$W$. Furthermore, $W_J$ is not
conjugate to any other visible subgroup of~$W$.
\end{lem}

Another important tool in our analysis of automorphisms is the
classification of involutions in Coxeter groups, due to
Richardson~\cite{RR}.

\begin{prop}\label{invclass}
Suppose that $w\in W$ is an involution. Then there is a $t\in W$ and
a spherical $I\subseteq \Pi$ such that $w=t^{-1}w_It$ with
$\ell(w)=\ell(w_I)+2\ell(t)$, and $w_I$ is central in $W_I$.
\end{prop}

\begin{proof}
See \cite[Proposition 5]{FH2}.
\end{proof}

\begin{defn}\label{-1type} {\rm We say that $I\subseteq\Pi$ is
\textit{of ($-$1)-type\/} if $W_I$ is finite and $w_I$ is central
in~$W_I$.}
\end{defn}

The reason for the terminology is that $I$ is of $(-1)$-type if and
only if there is an element of $W_I$ that acts on $V_I$ as
multiplication by~$-1$.

We need the following lemma.

\begin{lem}\label{parabconj}
Suppose that $I,\,J\subset\Pi$ with $I$ of ($-$1)-type, and suppose
that $t\in W$ has the property that $tw_It^{-1}\in W_J$. Then
$tW_It^{-1}\subseteq W_J$.
\end{lem}

\begin{proof} Let $a\in I$. Then $w_I(a)=-a$, and so
$(tw_It^{-1})(ta)=-ta$, whence it follows that either $ta$ or $-ta$ is
in the set $N(tw_It^{-1})$. But $N(tw_It^{-1})\subseteq\Phi_J$; so
$ta\in\Phi_J$, and therefore $tr_at^{-1}=r_{ta}\in W_J$. Since $W_I$
is generated by $\{\,r_a\mid a\in I\,\}$, the result follows.
\end{proof}

In particular, it follows from Lemma \ref{parabconj} that if $I,\,J$
are both of $(-1)$-type and $tw_It^{-1}=w_J$ then $tW_It^{-1}=W_J$.
Conversely, suppose that $tW_It^{-1}=W_J$, so that in fact
$dW_Id^{-1}=W_J$ for all $d$ in $W_JtW_I$ (which equals $tW_I$).
Taking $d$ to be the shortest element in $tW_I$,
Lemma~\ref{kilmoyerlem} yields that $dI=J$, and hence
$x\mapsto dxd^{-1}$ is a length-preserving isomorphism $W_I\to W_J$;
consequently $dw_Id^{-1}=w_J$. If $w_I,\,w_J$ are central in
$W_I,\,W_J$ we deduce that $tw_It^{-1}=w_J$. So we have proved the
following result.

\begin{lem}\label{WItoWJ} Suppose that $I,\,J$ are subsets of $\Pi$
that are both of ($-$1)-type. Then $\{\,t\in W\mid tw_It^{-1}=w_J\,\}=
\{\,t\in W\mid tW_It^{-1}=W_J\,\}$.
\end{lem}

\begin{prop}\label{Z(w)normal} Let $I\subset\Pi$ be of ($-$1)-type.
Then $W_I\subseteq \FC(w_I)$.
\end{prop}

\begin{proof} Let $F$ be a maximal finite subgroup of $W$ such that
$w_I\in H$. By Lemma~\ref{lemtits} there exist $t\in W$ and
$J\subseteq\Pi$ such that $tFt^{-1}=W_J$. By Lemma~\ref{parabconj}
and the fact that $w_I\in H$ it follows that
$tW_It^{-1}\subseteq W_J$. Hence $W_I\subseteq t^{-1}W_Jt=F$.
\end{proof}

\begin{prop}\label{key} Let $W,\,W'$ be Coxeter groups of finite
rank and $\alpha\colon W\to W'$ an isomorphism. Let $\Pi$ be the
set of simple roots corresponding to the distinguished generating
set of~$W$, and let $a\in\Pi$. If $r_a^\alpha$ is not a reflection
in $W'$ then the intersection of all maximal finite subgroups of
$W$ containing $r_a$ is a parabolic subgroup of order greater than~2.
\end{prop}

\begin{proof} Write $\Pi'$ for the set
of simple roots of~$W'$. Observe that Lemma~\ref{lemtits} and
Corollary~\ref{intersectionofparasispara} trivially imply that $\FC(r_a)$
is a parabolic subgroup of~$W$.

Since $r_a^\alpha$ is not a reflection it
follows from Proposition~\ref{invclass} that $r_a^\alpha=tw_It^{-1}$
for some $t\in W'$ and $I\subseteq\Pi'$ of $(-1)$-type and of rank at
least~2. Clearly $\FC(r_a)^\alpha=t\FC(w_I)t^{-1}$, and by
Proposition~\ref{Z(w)normal} we know that $W_I\subseteq \FC(w_I)$. Hence
$(tW_It^{-1})^{\alpha^{-1}}\subseteq \FC(r_a)$, so that $\FC(r_a)$ has
order greater than~2, as required.
\end{proof}

\begin{cor} \label{FCcrit1}
Let $w \in W$ be an involution
such that $\FC(w) = \langle w \rangle$
 and let $S \subseteq W$ be such that $(W,S)$ is a Coxeter system.
Then $w \in S^W$.
\end{cor}

\section{The finite continuation of a reflection}

Let $r$ be a reflection in~$W$.
Replacing $r$ by $wrw^{-1}$ replaces $\FC(r)$ by $w\FC(r)w^{-1}$, and so
choosing $w$ suitably enables us to assume that $\FC(r)=W_J$, a visible
parabolic subgroup. Furthermore, replacing $r$ by $trt^{-1}$ for suitable
$t\in W_J$ enables us to assume that $r=r_a$ for some $a\in J$.
(Note that these observations yield the first assertion of
Theorem~\ref{main}.)

\begin{prop}\label{finclos} Let $a\in J\subseteq\Pi$, and
suppose that $W_J$ is the intersection of all maximal finite
subgroups of $W$ containing~$r_a$. Then
$\{\,w\in W\mid wr_aw^{-1}\in W_J\,\}$ is a subset of the normalizer
of $W_J$ in~$W$. Thus each $W\!$-conjugate of $r_a$ in $W_J$ is
$N_W(W_J)$-conjugate to~$r_a$, and $C_W(r_a)\subseteq N_W(W_J)$.
Moreover, if $b\in\Pi\setminus J$ is such that $W_{J\cup\{b\}}$ is
infinite then $m_{bc}=\infty$ for all $c\in J$ such that $r_c$ is
conjugate to $r_a$ in~$W$.
\end{prop}

\begin{proof}
Let $\mathscr{S}$ be the set of all maximal finite subgroups of $W$
containing~$r_a$, so that $W_J=\FC(r_a)=\bigcap_{F\in\mathscr{S}}F$.
Suppose that $w\in W$ satisfies $wr_aw^{-1}\in W_J$, and let
$F\in\mathscr{S}$. Then $wr_aw^{-1}\in W_J\subseteq F$, and so
$r_a\in w^{-1}Fw$. Thus $w^{-1}Fw$ is a maximal finite subgroup of
$W$ containing~$r_a$, whence $w^{-1}Fw\in\mathscr{S}$. So
$$
\bigcap_{F\in\mathscr{S}}F\subseteq\bigcap_{F\in\mathscr{S}}w^{-1}Fw
$$
and so $W_J\subseteq w^{-1}W_Jw$. Since $W_J$ is finite it follows
that $w\in N_W(W_J)$.

Suppose that $c\in J$ with $r_c=wr_aw^{-1}$ for some $w\in W$.
Clearly $F\mapsto wFw^{-1}$ is a bijection from the set of maximal
finite subgroups of $W$ containing $r_a$ to the set of maximal
finite subgroups of $W$ containing~$r_c$, and so
$\FC(r_c)=w\FC(r_a)w^{-1}$. But $w\FC(r_a)w^{-1}=wW_Jw^{-1}=W_J$ by
the first part of the proof, and so $\FC(r_c)=W_J$. Now suppose that
$b\in\Pi\setminus J$ with $m_{cb}<\infty$. Then $W_{\{c,b\}}$ is
finite, and so contained in a maximal finite subgroup~$F$. Since
$r_c\in F$ we must have $\FC(r_c)\subseteq F$. It follows that the
finite group $F$ contains both $W_J$ and~$r_b$, and therefore
$W_{J\cup\{b\}}$ is finite.
\end{proof}

Assume, as in Proposition \ref{finclos}, that
$a\in J\subseteq\Pi$ and $W_J=\FC(r_a)$, and suppose now that
$J\ne\{a\}$. Suppose that $L\subseteq\Pi$ is such that $J\subseteq L$
and $W_L$ is finite. Then
$W_L$ is a finite Coxeter group possessing a visible parabolic
subgroup $W_J$ of rank greater than~1 that is normalized by the
centralizer of some simple reflection $r_a\in W_J$. Indeed, $W_J$ is
normalized by all $w\in W_L$ such that $wr_aw^{-1}\in W_J$.
Equivalently, by Lemma~\ref{somebasicfacts}~(3),
$\{\,w\in W_L\mid wa\in\Phi_J\,\}\subseteq N_W(W_J)$.
This is a very restrictive condition, which we now proceed to examine
with a case-by-case investigation of the different types of finite
Coxeter groups. For the course of this investigation, we can (and
shall) assume that $L=\Pi$.

So we assume for now that $W$ is a finite Coxeter group of rank~$n$,
and our aim is to find all examples of the following phenomenon:
there exist $\{a\}\subsetneqq J\subseteq\Pi$ such that the set
$Q=\{\,w\in W\mid wa\in\Phi_J\,\}$ is a subset of $N_W(W_J)$.
We assume that $J\ne\Pi$, since the condition is trivially
satisfied otherwise.

If $K\subseteq\Pi$ is a component of the Coxeter graph
such that $J\cap K=\emptyset$ then $W_K$ is a direct factor of
$N_W(W_J)$; moreover, $Q=(Q\cap W_{\Pi\setminus K})W_K$. So
removing $K$ from the graph will have no bearing on whether or
not the condition $Q\subseteq N_W(W_J)$ holds. So we assume that
there are no such components of~$\Pi$. Exactly the same comments
apply for a component $K$ of $\Pi$ such that $K\subseteq J$.
So we also assume that there are none of these.

Assume that $\{a\}\subsetneqq J\subsetneqq\Pi$ and
$Q\subseteq N_W(W_J)$. Suppose
that $K\subseteq\Pi$ is a component of the Coxeter graph
such that $a\notin K$. Then $r_ba=a$ for all $b\in K$; so
$r_b\in Q\subseteq N_W(W_J)$, and it follows that $r_bc\in\Phi_J$
whenever $c\in J$. If $b\cdot c\ne 0$ then $b$ is in the support
of~$r_bc$, and so $r_bc\in\Phi_J$ implies $b\in J$. Since $K$
is connected it follows that if $K$ contains any element of $J$ then
$K\subseteq J$. So either $K\cap J=\emptyset$ or $K\subseteq J$.
But we have assumed that there are no such components. So the
component of $\Pi$ that contains $a$ is the only component;
that is, $\Pi$ is irreducible.

Observe that the group $\Stab(a)=\{\,w\in W\mid wa=a\,\}$ is a subset
of~$Q$ and hence of $N_W(W_J)$. Note also that
$N_W(W_J)=\{\,w\in W\mid w\Phi_J=\Phi_J\,\}$, which is also the
stabilizer of the subspace~$V_J$ (since $V_J$ is the subspace
spanned by $\Phi_J$ and $\Phi_J=V_J\cap\Phi$). Now
$\Stab(a)$ is a parabolic subgroup of $W$ whose root system is
$\Phi\cap a^\perp$, and the following table gives the type of this
root system in all cases.
$$
\vbox{\offinterlineskip
\halign{\quad$\hfil#\hfil$\quad&\vrule#&\quad$\hfil#\hfil$\cr
W& depth 4.5 pt&\Stab(a)\cr
\noalign{\hrule}
A_n&height 9.5 pt depth 3.5 pt&A_{n-2}\cr
C_n&height 8.5 pt depth 3.5 pt&C_{n-2}+A_1\cr
C_n&height 8.5 pt depth 3.5 pt&C_{n-1}\cr
D_n&height 8.5 pt depth 3.5 pt&D_{n-2}+A_1\cr
F_4&height 8.5 pt depth 3.5 pt&C_3\cr
E_6&height 8.5 pt depth 3.5 pt&A_5\cr}}
\qquad\qquad\qquad
\vbox{\offinterlineskip
\halign{$\hfil#\hfil$\quad&\vrule#&\quad$\hfil#\hfil$\cr
W& depth 4.5 pt&\Stab(a)\cr
\noalign{\hrule}
E_7&height 8.5 pt depth 3.5 pt&D_6\cr
E_8&height 8.5 pt depth 3.5 pt&E_7\cr
H_3&height 8.5 pt depth 3.5 pt&A_1+A_1\cr
H_4&height 8.5 pt depth 3.5 pt&H_3\cr
I_2(2k)&height 8.5 pt depth 3.5 pt&A_1\cr
I_2(2k+1)&height 8.5 pt depth 3.5 pt&\emptyset\cr}}
$$
(For $C_n$ there are two $W$ orbits of roots, giving two possibilities
for $\Stab(a)$. For $F_4$ and $I_2(2k)$ there are also two $W$-orbits
of roots, but $\Stab(a)$ has the same type of root system whichever
orbit $a$ belongs to.) Since each irreducible constituent of its root
system spans an irreducible $\Stab(a)$-submodule of~$V$, the table
shows that as a $\Stab(a)$-module, $V$ has composition length two or
three or (in one case only) four: $a$ itself spans a trivial $\Stab(a)$
submodule of dimension~1,
and $a^\perp$ is either irreducible of dimension~$n-1$ (for types
$F_4$, $E_6$, $E_7$, $E_8$, $H_4$, $I_2(2k)$ and one of the $C_n$
possibilities), or the direct sum of irreducibles of dimensions
1~and~$n-2$ (for types $A_n$, $C_n$, $D_n$ when $n>4$, $H_3$ and
$I_2(2k+1)$), or the direct sum of three irreducibles of dimension~1
(for type~$D_4$). Furthermore, the summands of $a^\perp$ are
pairwise nonisomorphic as $\Stab(a)$ modules, since even if they are
of the same type their centralizers in $\Stab(a)$ are different.

Since $\{a\}\subsetneqq J\subsetneqq\Pi$ and $V_J$ is
$\Stab(a)$-invariant, we see that
$a^\perp=(V_J\cap a^\perp)\oplus V_J^\perp$, with both summands
nonzero $\Stab(a)$-modules. So $\Pi$ is of type $A_n$, $C_n$, $D_n$
or $H_3$. Furthermore, except in type $D_4$, the two direct summands
of $a^\perp$ are irreducible and not isomorphic, and are therefore
the only proper $\Stab(a)$-submodules of~$a^\perp$. We conclude that
$V_J$ is spanned by $a$ and one of the summands of~$a^\perp$, while
$V_J^\perp$ is the other summand. In type $D_4$ we similarly deduce
that $V_J$ is spanned by $a$ and one or two of the three 1-dimensional
summands of~$a^\perp$, and, correspondingly, $V_J^\perp$ is of either
of type $A_1+A_1$ or of type~$A_1$.

If $\Pi$ is of type $A_n$ then one of the summands of $a^\perp$ is of
type $A_{n-2}$ while the other is a trivial 1-dimensional
$\Stab(a)$-module. If $V_J^\perp$ is of type $A_{n-2}$ then $V_J$ must
be of type $A_1$, since the orthogonal complement of a subsystem of
type $A_{n-2}$ in $A_n$ contains only a rank~1 root system. This
contradicts the assumption that $\{a\}\subsetneqq J$. So $J$ is
of type $A_1+A_{n-2}$. Since $W_J$ is visible, we deduce that $a$ is
an end node of the $A_n$ diagram, and the node adjacent to $a$ is the
unique simple root not in~$J$. However, if $n>3$ then the maximal
length element of $W$ is in~$Q$ but not in the normalizer of~$W_J$.
So $n=3$ and $J=\{a,c\}$, where $c$ is the other end node. It is
readily checked that $Q$ has 8 elements and coincides with $N_W(W_J)$
(which is generated by $W_J$ and an element that interchanges
$a$~and~$c$).

If $\Pi$ is of type $C_n$ then one summand of $a^\perp$ is of type
$C_{n-2}$ and the other of type~$A_1$. The roots in the $A_1$ summand
are in the same $W$-orbit as~$a$. If $V_J^\perp$ is the $A_1$
component of $a^\perp$ then $V_J=(V_J^\perp)^\perp$ is of type
$C_{n-2}+A_1$. This determines $J$ uniquely, since $W_J$ is visible.
If $n\ge 4$ and $w$ is the longest element of the visible parabolic
subgroup of type $A_{n-1}$, then $-wa\in\Pi\setminus\{b\}=J$, but
$w\notin N_W(W_J)$. This contradicts the fact that
$Q\subseteq N_W(W_J)$. So $n=3$, and the elements of $J$ are the end
nodes $a,\,c$ of the $C_3$ diagram, the middle node $b$ being in the
same $W$-orbit as~$a$. Since $\Phi_J=\{\pm a,\pm c\}$ and $c$ is not
in the same $W$-orbit as $a$~and~$b$ we deduce that
$Q=\{\,w\in W\mid wa=\pm a\,\}$. Furthermore, of the 6 roots in
the $W$-orbit of~$c$, only $c$ and $-c$ are orthogonal to~$a$. So
if $wa=\pm a$ then $wc=\pm c$. Thus if $w\in Q$ then
$w\Phi_J=\Phi_J$, as required.

Continuing the discussion of~$C_n$, suppose now that $V_J^\perp$ is
the $C_{n-2}$ component of~$a^\perp$. Then $V_J=(V_J^\perp)^\perp$ is
of type~$C_2$. Writing $J=\{a,b\}$, the fact that $\Stab(a)$ is of
type $A_1+C_{n-2}$ means that it is $b$ rather than $a$ that is
the end node of the $C_n$ diagram. If we put $c=r_ba$ then $\{\pm c\}$
is the component of $\Phi\cap a^\perp$ of type~$A_1$. It follows that
$\{\pm a\}=\{\pm r_bc\}$ is the $A_1$-component of
$\Phi\cap(r_ba)^\perp=\Phi\cap c^\perp$. We see that
$\Stab(a)=\langle r_c\rangle\times W'$ and
$\Stab(c)=\langle r_a\rangle\times W'$, where $W'$ is a parabolic
(not visible) subgroup of~$W$ of type~$C_{n-2}$. Indeed, the root
system of $W'$ is $\Phi\cap V_J^\perp$. The roots in $\Phi_J$ that
are in the same $W$-orbit as $a$ are $\pm a$ and $\pm c$, and so
$$
Q=\{1,r_a,r_b,r_br_a\}\Stab(a)=\{1,r_a,r_b,r_br_a\}\{1,r_c\}W'
=W_JW'.
$$
Hence  our requirement that $Q$ stabilizes $\Phi_J=\{\pm a,\pm b\}$
is indeed satisfied.

If $\Pi$ is of type $D_n$ with $n>4$ then one summand of $a^\perp$
is of type $D_{n-2}$ and the other of type~$A_1$. The roots orthogonal
to a $D_{n-2}$ subsystem form a system of type~$A_1+A_1$. There
are in fact two $W$-orbits of parabolic $A_1+A_1$ subsystems,
and the orthogonal complement of a $D_{n-2}$ is perhaps better thought
of as type~$D_2$, since the visible parabolic in this orbit
corresponds to the two nodes of the diagram that form the fork. So
if $V_J^\perp$ is the $D_{n-2}$ summand of $a^\perp$ then $J=\{a,b\}$
consists to the two nodes of valency~1 that are adjacent to~$c$, the
node of valency~3. A similar statement applies for $D_4$ in the case
that $V_J^\perp$ is of type $A_1+A_1$. In both cases the element
$w=r_cr_ar_br_c\in W$ satisfies $wa=b$ and $wb=a$, and since
$\Phi_J=\{\pm a,\pm b\}$ we see that
$Q=\{1,r_a,w,wr_a\}\Stab(a)$.
But $\Stab(a)=\langle r_b\rangle\times W'$ and
$\Stab(b)=\langle r_a\rangle\times W'$, where $W'$ is the parabolic
subgroup corresponding to the subspace~$V_J^\perp$, and it follows
readily that $Q$ stabilizes $\Phi_J=\{\pm a,\pm b\}$, as required.

Continuing the discussion of~$D_n$, where $n\ge 4$, suppose now that
$V_J^\perp$ is an $A_1$ component of~$a^\perp$. Then
$V_J=(V_J^\perp)^\perp$ is of type $A_1+D_{n-2}$. But the maximal
length element of a visible $A_{n-1}$ subsystem containing~$a$ takes
$a$ to an element of $\Phi_J$ but does not normalize~$W_J$. So our
requirement that $Q\subseteq N_W(W_J)$ is not met.

Finally, suppose that $\Pi$ is of type $H_3$, so that $\Stab(a)$ is
of type $A_1+A_1$. Then $V_J^\perp$ is of type~$A_1$, and hence $J$
is of type $A_1+A_1$. Let $J=\{a,c\}$, and note that $c=wa$ for
some~$w\in W$. Since $N_W(W_J)$ is generated by $W_J$ and the central
involution of~$W$, we see that $c$ is not in the $N_W(W_J)$-orbit
of~$a$. Hence the element $w$ above is in $Q$ but not in $N_W(W_J)$,
and so our requirements are not met.

We have thus established the following result.

\begin{prop}\label{cases}
Let $\Pi$ be the set of simple roots for the finite irreducible
Coxeter group~$W$, and suppose that $a\in J\subseteq\Pi$.
Then $\{\,w\in W\mid wa\in\Phi_J\,\}$ is a subset of $N_W(W_J)$ if and
only if one of the following situations occurs\textup{:}
\begin{enumerate}
\item $J=\{a\}$\textup{;}
\item $J=\Pi$\textup{;}
\item $\Pi=\{a,b,c\}$ is of type $C_3$, with $m_{ab}=3$ and $m_{bc}=4$,
and $J=\{a,c\}$ of type $A_1+A_1$\textup{;}
\item $\Pi$ is of type~$D_n$ or $A_3$, and $J=\{a,b\}$, where $a$ and
$b$ are end nodes that are both adjacent to some $c\in\Pi$\textup{;}
\item $\Pi$ is of type~$C_n$ and $J=\{a,b\}$ is of type $C_2$,
with $b$ an end node of~$\Pi$.
\end{enumerate}
\end{prop}

We return now to investigation of an arbitrary finite rank Coxeter
group~$W$. The next proposition is an immediate consequence of
Proposition~\ref{cases} and the discussion preceding it.

\begin{prop}\label{finclos2}
Let $a\in J\subseteq L\subseteq\Pi$, and suppose that the group
$W_L$ is finite and that $W_J = \FC(r_a)$. 
Let $J_0$ be the
component of $J$ containing~$a$ and $L_0$ the component of
$L$ containing~$J_0$. Then every component of $J$ that
is not contained in $L_0$ is a component of~$L$.
Furthermore, if $\{a\}\ne J\cap L_0\ne L_0$ then $J\cap L_0=\{a,b\}$
for some~$b$, and one of the following alternatives occurs\/\textup{:}
\begin{enumerate}
\item $L_0=\{a,c,b\}$ is of type~$C_3$, with $m_{ac}=3$ and
$m_{cb}=4$\textup{;}
\item $L_0$ is of type~$C_n$ for some~$n\ge 3$, with $b$ an end node
and $J_0=\{a,b\}$ of type~$C_2$\textup{;}
\item $L_0$ is of type $A_3$ or type $D_n$ for some $n\ge 4$, the
nodes $a$~and~$b$ having valency~1 and sharing a common neighbour.
\end{enumerate}
\end{prop}

One of the ingredients of alternative~(2) of Proposition~\ref{finclos2}
is that the component of $\FC(r_a)$ containing~$a$ is of type~$C_2$.
We shall see that when this situation arises, $\Odd(a)$
has a focus in~$\Pi$.

\begin{prop}\label{alt2}Suppose that $a\in J\subseteq\Pi$ with $W_J=\FC(r_a)$,
and let $J_0$ be the component of $J$ containing~$a$. Suppose that
$J_0=\{a,b\}$ is of type~$C_2$.
Then either $\Odd(a)\cup\{b\}$ is a spherical component of $\EOdd(a)$,
or else $(a,b)$ is a focus of $\Odd(a)$~in~$\Pi$.
\end{prop}

\begin{proof}
We use induction on $k$ to prove that for all $k\ge 2$, if
$b=c_1,\,a=c_2,\,c_3,\,\ldots,\,c_k$ are simple roots satisfying
\begin{enumerate}
\item
$2<m_{c_ic_{i+1}}<\infty$ for all $i\in\{1,2,\ldots,k-1\}$, and
\item
$c_1,\,c_2,\,\ldots,\,c_k$ are distinct from each other,
\end{enumerate}
then $\{c_1,c_2,\ldots,c_k\}$ forms a system of type~$C_k$. The case $k=2$ is
immediately true.

Suppose that $k>2$. The inductive hypothesis tells us that
$\{c_1,c_2,\ldots,c_{k-1}\}$ is of type~$C_{k-1}$. The element
$w=v[c_{k-1},\{c_{k-2}\}]\cdots v[c_4,\{c_3\}]v[c_3,\{c_2\}]$ has the property
that $wa=wc_2=c_{k-1}$, and so if we write $d=wb$ then
\[
r_d=wr_bw^{-1}\in w\FC(r_a)w^{-1}=\FC(c_{k-1}),
\]
since it is given that $b\in\FC(r_a)$. But $W_{\{c_{k-1},c_k\}}$ is finite, and
so it follows that $\{r_d,r_{c_{k-1}},r_{c_k}\}$ generates a finite group. Now
$d\cdot c_{k-1}=b\cdot a=-\cos(\pi/4)$ and $c_{k-1}\cdot c_k=-\cos(\pi/m)$
for some $m>2$. If $m\ge 4$ then
\[
c_k\cdot d=c_k\cdot\Bigl(b+\sqrt 2\sum_{i=2}^{k-1}c_i\Bigr)
\le -\sqrt 2(c_k\cdot c_{k-1})\le-1,
\]
whence the reflection subgroup $W_{\{r_d,r_{c_k}\}}$ is infinite (by
Theorem~\ref{base}), a contradiction. So $m=3$. If $m_{c_ic_k}>2$ for any
$i\in\{1,2,\ldots,k-2\}$ then
\[
c_k\cdot d\le \sqrt 2(c_k\cdot c_{k-1})+c_k\cdot c_i<-1,
\]
again giving a contradiction. So $m_{c_kc_i}=2$ for all
$i\in\{1,2,\ldots,k-2\}$ and $m_{c_kc_{k-1}}=3$, and since
$\{c_1,c_2,\ldots,c_{k-1}\}$ is a system of type~$C_{k-1}$ it follows that
$\{c_1,c_2,\ldots,c_k\}$ is a system of type~$C_k$, as claimed.

If there were $c,\,d\in\Odd(a)$ with $3<m_{cd}<\infty$ then $b$ together
with a minimal length odd-labelled path from $a$ to $\{c,d\}$ would
yield $c_1,\,c_2,\,\ldots,\,c_k\in\Pi$ satisfying (1) and (2) above and not
forming a system of type~$C_k$, contradicting the result proved above.
The same argument yields a contradiction if $c\in\Odd(a)$ and
$d\in\Pi\setminus\Odd(a)$ with $3<m_{cd}<\infty$, unless $\{c,d\}=\{a,b\}$.
So all edge labels in $\Odd(a)$ are~$3$, if $c,\,d\in\Odd(a)$ are not adjacent
in $\Odd(a)$ then $m_{cd}\in\{2,\infty\}$, and if $c\in\Odd(a)$
and $d\in\Pi\setminus\Odd(a)$ then $m_{cd}\in\{2,\infty\}$ unless
$\{c,d\}=\{a,b\}$. Furthermore, any circuit in $\Odd(a)$ would similarly
yield a contradiction (by combining the circuit with a minimal
finite-labelled path connecting it to~$b$). So $\Odd(a)$ is tree.

For each $c\in\Odd(a)$ let $C[b..c]\subseteq\Pi$ consist of $b$ and the unique
path from $a$ to $c$ in $\Odd(a)$. The discussion above shows that $C[b..c]$ is
always of type~$C$. Now suppose that $c\in\Odd(a)$ and $e\in\Pi\setminus C[b..c]$
with $m_{ce}=2$. Write $C[b..c]=\{c_1,c_2,\cdots,c_k\}$, with $c_1=b$ and
$c_k=c$, and let $d=b+\sqrt 2\sum_{i=2}^kc_i$. An argument similar to one used
above shows that $r_d\in\FC(c)$, and hence $W_{\{d,c,e\}}$ is finite.
So $d\cdot e>-1$. If $c_i\cdot e\ne0$ then $c_i\cdot e\le-1/2$; so it follows
that there is at most one~$i$ with $c_i\cdot e\ne0$. Suppose, for a contradiction,
that there is exactly one such~$i$. If $i>1$ then $d\cdot e=\sqrt 2(c_i\cdot e)$,
and so $c_i\cdot e>-1/\sqrt 2$. Hence $m_{c_ie}=3$, and $d\cdot e=-1/\sqrt 2$.
But this means that the edges $\{c,d\}$ and $\{d,e\}$ of the Coxeter graph
of $\{d,c,e\}$ are both labelled~$4$, contradicting the fact that
$W_{\{d,c,e\}}$ is finite. So we must have $i=1$, and finiteness of $W_{\{c,d,e\}}$
forces $b\cdot e=d\cdot e=-1/2$. But now if we put $L=\{e\}\cup J$ then, in the
notation of Proposition~\ref{finclos2}, we have that $L_0=\{e,b,a\}$ is of type
$C_3$ with $J\cap L_0=\{b,a\}$ of type~$C_2$, and Proposition~\ref{finclos2} shows
that this is not possible. We conclude that if $e\in\Pi$ has the property that
$m_{ce}=2$ for some $c\in\Odd(a)$ then $m_{de}=2$ for all~$d\in C[b..c]$. In
particular, if $e\in\Pi\setminus(\Odd(a)\cup\{b\})$ and $m_{ce}\ne\infty$ for some
$c\in\Odd(a)$ then $m_{ce}=2$, as shown above, and so
$m_{de}=2$ for all $d\in C[b..c]$.

All that remains to prove now is that if $c,\,d\in\Odd(a)$ with $c\notin C[b..d]$
and $d\notin C[b..c]$, then $m_{cd}=\infty$. Since $c$ and $d$ are not adjacent in
$\Odd(a)$ the only alternative is that $m_{cd}=2$; so suppose, for a contradiction,
that this holds. Choose the vertex $e\in\Odd(a)$ on the (unique) path from
$c$ to~$d$ such that the distance from $e$ to~$a$ is minimal. Let $c',\,d'$ be the
neighbours of $e$ in the path from $c$~to~$d$, with $c'$ between $e$~and~$c$ and
$d'$ between $e$~and~$d$. Then $c'\in C[b..c]$, and since $m_{cd}=2$ it follows
that $m_{c'd}=2$. Now since $d'\in C[b..d]$ and $m_{c'd}=2$ it follows that
$m_{c'd'}=2$. Thus the set $L\subseteq\Pi$ consisting of $c'$ and $d'$ and the
vertices on the path from $a$~to~$e$ form a system of type $D$ (or $A_3$ if~$e=a$).
So $L$ is spherical, and since $b\in\FC(r_a)$ it follows that $L\cup\{b\}$ is
spherical also. But this is impossible since $L\cup\{b\}$ is connected, has an edge
labelled~4 (namely, $\{b,a\}$), and has a vertex of valency~3 (namely~$e$).
\end{proof}

The situation of alternative~(3) of Proposition~\ref{finclos2} is very similar
to that of alternative~(2), and in this case it turns out that $\Odd(a)$ has a
half-focus in~$\Pi$.

\begin{prop}\label{alt3}Suppose that $a\in J\subseteq\Pi$ with $W_J=\FC(r_a)$
and $\{a\}$ a component of~$J$, and suppose that $J\cap\Odd(a)\ne\{a\}$. Then
either $\Odd(a)$ is a spherical component of~$\EOdd(a)$, or else
there exists an element $b\in\Odd(a)$ such that $\{a,b\}$ is a half focus of $\Odd(a)$
in~$\Pi$.
\end{prop}

\begin{proof}
Let $b\in(J\cap \Odd(a))\setminus\{a\}$, and let $w\in W$ with $wa=b$. Then
$w\in N_W(W_J)$, by Proposition~\ref{finclos}, and so
\[
\FC(r_b)=\FC(wr_aw^{-1})=w\FC(r_a)w^{-1}=wW_Jw^{-1}=W_J.
\]
Moreover, $w\Phi_J=\Phi_J$, and since $a\cdot c=0$ for all
$c\in\Phi_J\setminus\{a\}$, it follows that $wa\cdot d=0$ for all
$d\in\Phi_J\setminus\{wa\}$. So $\{b\}$ is a component of~$J$. Note that
$m_{ab}=2$, since $a$~and~$b$ are in different components of~$J$.

Let $c\in\Pi\setminus\{a,b\}$, and suppose first of all that $2<m_{bc}<\infty$.
Since $\{b,c\}$ is spherical and $W_J=\FC(r_b)$ it follows that $J\cup\{c\}$ is
spherical. Let $L=J\cup\{c\}$ and let $L_0$ be the component of $L$ containing~$a$.
By Proposition~\ref{finclos2}, every component of $J$ that is not contained in
$L_0$ is a component of~$L$. But $b$ is adjacent to $c$~in~$L$; so $\{b\}$ is not
a component of~$L$, and it follows that~$b\in L_0$. Now $\{a\}\ne J\cap L_0$,
since $b\in J\cap L_0$, and $J\cap L_0\ne L_0$, since $c\in L_0$ and $c\notin J$
(since $\{b\}$ is a component of~$J$). Furthermore, the conditions of
alternative~(2) of Proposition~\ref{finclos2} are not satisfied, since $a$~and~$b$
are not adjacent in~$J$. So either alternative~(1) or alternative~(3) must hold,
and since $c$ is the only element of~$L$ not in~$J$ it follows that
$L_0=\{a,c,b\}$, with $m_{ac}=3$. But a symmetrical argument, with the roles
of $a$~and~$b$ reversed, shows that every $d\in\Pi$ with $2<m_{ad}<\infty$ has the
property that $m_{bd}=3$. So $m_{ac}=m_{bc}=3$, and $\{a,c,b\}$ is of type~$A_3$.

Now suppose that $m_{bc}=2$. Again since $\{b,c\}$ is spherical it follows that
$J\cup\{c\}$ is spherical, and so $m_{ac}<\infty$. If $m_{ac}>2$ then, as we have
just observed, it follows that $m_{bc}=3$, contrary to our assumption that
$m_{bc}=2$. So $m_{ac}=m_{bc}=2$, and we have now shown that whenever
$m_{bc}<\infty$ we have $m_{ac}=m_{bc}\in\{2,3\}$. Since a symmetrical argument
gives the same conclusion whenever $m_{ac}<\infty$, we conclude also that
$m_{ac}=\infty$ if and only if $m_{bc}=\infty$.

We now use induction on $k$ to prove that for all $k\ge 3$, if
$b=c_1,\,a=c_2,\,c_3,\,\ldots,\,c_k$ are simple roots satisfying
\begin{enumerate}
\item
$2<m_{c_ic_{i+1}}<\infty$ for all $i\in\{2,3,\ldots,k-1\}$, and
\item
$c_1,\,c_2,\,\ldots,\,c_k$ are distinct from each other,
\end{enumerate}
then $\{c_1,c_2,\ldots,c_k\}$ forms a system of type~$D_k$ or~$A_3$. The case
$k=3$ follows from what we have proved above.

Suppose that $k>3$. The inductive hypothesis tells us that
$\{c_1,c_2,\ldots,c_{k-1}\}$ is of type~$D_{k-1}$ (or $A_3$ if $k=4$). The
element $w=v[c_{k-1},\{c_{k-2}\}]\cdots v[c_4,\{c_3\}]v[c_3,\{c_2\}]$ has the property
that $wa=wc_2=c_{k-1}$, and so if we write $d=wb$ then
\[
r_d=wr_{c_1}w^{-1}\in w\FC(r_a)w^{-1}=\FC(c_{k-1}),
\]
since it is given that $b\in\FC(r_a)$. But $W_{\{c_{k-1},c_k\}}$ is finite, and so
it follows that $\{r_d,r_{c_{k-1}},r_{c_k}\}$ generates a finite group. Now
$d\cdot c_{k-1}=b\cdot a=0$ and $c_{k-1}\cdot c_k=-\cos(\pi/m)$
for some $m>2$. If $c_k\cdot c_i\ne 0$ for some $i\in\{1,2,\ldots,k-2\}$ then
\[
c_k\cdot d=c_k\cdot\Bigl(c_1+c_2+c_{k-1}+2\sum_{i=3}^{k-2}c_i\Bigr)\le
-(\tfrac12+\cos\tfrac\pi m)\le-1,
\]
whence the reflection subgroup $W_{\{r_d,r_{c_k}\}}$ is infinite (by
Theorem~\ref{base}), a contradiction. So
$c_k\cdot d=c_{k-1}\cdot c_k=-\cos\frac\pi m$. Since the reflection subgroup
generated by $\{r_d,r_{c_{k-1}},r_{c_k}\}$ is finite it follows that $m=3$.
So we have shown that $m_{c_kc_{k-1}}=3$ and $m_{c_kc_i}=2$ for $i<k-1$,
and since $\{c_1,c_2,\ldots,c_{k-1}\}$ is a system of type~$D_{k-1}$ it follows
that $\{c_1,c_2,\ldots,c_k\}$ is a system of type~$D_k$, as claimed.

Note that $\Odd(a)\setminus\{b\}$ and $\Odd(a)\setminus\{a\}$ are both
connected, since each element $c\in\Odd(a)$ that is is adjacent to~$a$ is
also adjacent to~$b$, and vice versa.
If there were $c,\,d\in\Odd(a)\setminus\{b\}$ with $3<m_{cd}<\infty$ then
$b$ together with a minimal length odd-labelled path from $a$ to $\{c,d\}$ would
yield $c_1,\,c_2,\,\ldots,\,c_k\in\Pi$ satisfying (1) and (2) above and not
forming a system of type~$D_k$, contradicting the result proved above.
The same argument yields a contradiction whenever $c\in\Odd(a)\setminus\{a,b\}$
and $d\in\Pi\setminus\Odd(a)$ with $3<m_{cd}<\infty$.
So all edge labels in $\Odd(a)$ are~$3$, if $c,\,d\in\Odd(a)$ are not adjacent
in $\Odd(a)$ then $m_{cd}\in\{2,\infty\}$, and if $c\in\Odd(a)$
and $d\in\Pi\setminus\Odd(a)$ then $m_{cd}\in\{2,\infty\}$. Furthermore, any
circuit in $\Odd(a)\setminus\{b\}$ would similarly
yield a contradiction (by combining the circuit with a minimal
finite-labelled path connecting it to~$b$). So $\Odd(a)\setminus\{b\}$ is tree.
Of course, $\Odd(b)\setminus\{a\}$ is also a tree, by the same argument.

For each $c\in\Odd(a)\setminus\{a,b\}$ let $D[a,b..c]\subseteq\Pi$ consist of
$b$ and the unique path from $a$ to $c$ in $\Odd(a)\setminus\{b\}$. The
discussion above shows that $D[a,b..c]$ is of type~$D$. Now suppose that
$c\in\Odd(a)\setminus\{a,b\}$ and $e\in\Pi\setminus D[a,b..c]$ with $m_{ce}=2$.
Write $D[a,b..c]=\{c_1,c_2,\cdots,c_k\}$, with $c_1=b$, $c_2=a$ and
$c_k=c$, and let $d=c_1+c_2+c_k+2\sum_{i=3}^{k-1}c_i$. An argument similar to one
used above shows that $r_d\in\FC(c)$, and hence $W_{\{d,c,e\}}$ is finite.
So $d\cdot e>-1$. If $c_i\cdot e\ne0$ then $c_i\cdot e\le-1/2$; so it follows
that $\{\,i\mid c_i\cdot e\ne0\,\}$ is a subset of $\{1,2,k\}$ with at most
one element. But $c_k\cdot e=0$ since $m_{ce}=2$, and $c_1\cdot e=c_2\cdot e$
since $m_{af}=m_{bf}$ for all~$f\in\Pi$. So $c_i\cdot e=0$ for all
$i\in\{1,2,\ldots,k\}$. In particular, if $e\in\Pi\setminus\Odd(a)$ and
$m_{ce}\ne\infty$ for some $c\in\Odd(a)$ then $m_{ce}=2$, as shown above,
and so $m_{de}=2$ for all $d\in D[a,b..c]$.

All that remains to prove now is that if $c,\,d\in\Odd(a)\setminus\{a,b\}$ with
$c\notin D[a,b..d]$ and $d\notin D[a,b..c]$, then $m_{cd}=\infty$. Since $c$ and
$d$ are not adjacent in $\Odd(a)$ the only alternative is that $m_{cd}=2$; so
suppose, for a contradiction, that this holds. Choose the vertex
$e\in\Odd(a)\setminus\{b\}$ on the (unique) path from $c$ to~$d$ such that the
distance from $e$ to~$a$ is minimal. Let $c',\,d'$ be the neighbours of $e$ in
the path from $c$~to~$d$, with $c'$ between $e$~and~$c$ and $d'$ between
$e$~and~$d$. Then $c'\in C[b..c]$, and since $m_{cd}=2$ it follows that
$m_{c'd}=2$. Now since $d'\in C[b..d]$ and $m_{c'd}=2$ it follows that $m_{c'd'}=2$.
Thus the set $L\subseteq\Pi$ consisting of $c'$ and $d'$ and the vertices on the
path from $a$~to~$e$ form a system of type $D_k$, or $A_3$ if~$e=a$.
So $L$ is spherical, and since $b\in\FC(r_a)$ it follows that $L\cup\{b\}$ is
spherical also. If $L=A_3$ then $m_{ac'}=m_{ad'}=3$, and since $m_{bc'}=m_{ac'}$
and $m_{bd'}=m_{ad'}$ we see that $L\cup\{b\}$ is of type $\widetilde A_3$,
contradicting the fact that $L\cup\{b\}$ is spherical. Similarly if
$L$ is of type~$D_k$ then $L\cup\{b\}$ is of type $\widetilde D_k$, again
giving a contradiction.
\end{proof}

We also need to obtain further information about the situation of
alternative~(1) of Proposition~\ref{finclos2}.
So for the next three lemmas we assume that $a\in J\subseteq L\subseteq\Pi$
with $L$ spherical and $W_J=\FC(r_a)$, and there exist $b\in J$ and
$c\in L\setminus J$ such that $L_0=\{a,c,b\}$ is a component of $L$ of
type~$C_3$, with $m_{ac}=3$ and $m_{cb}=4$.

\begin{lem}\label{alt1a}
For all $e\in\Pi\setminus\{a,c,b\}$, either
$m_{ce}=m_{ae}=m_{be}=2$ or $m_{ce}=\infty$.
Moreover, $J\cap\Odd(a)=\{a\}$.
\end{lem}

\begin{proof}
If $J\cap\Odd(a)\ne\{a\}$ then, since $\{a\}$ is a component of~$J$,
Proposition~\ref{alt3} applies, and it follows in particular that no
vertex in $\Odd(a)$ lies on an edge with finite label different from~3.
This contradicts $m_{bc}=4$. So $J\cap\Odd(a)=\{a\}$.

Suppose that $e\in\Pi\setminus\{a,c,b\}$ with $m_{ce}<\infty$. The group
$r_cr_aW_{\{c,e\}}r_ar_c$ is finite and contains $r_cr_aW_{\{c\}}r_ar_c=W_{\{a\}}$;
so there exists a maximal finite subgroup $G$ of~$W$ containing~$r_a$ and the
reflection along~$(r_cr_a)e$. Since $r_b\in\FC(r_a)\subseteq G$ it
follows that $W_{\{b,(r_cr_a)e\}}$ is finite, and hence so is
$W_{\{(r_ar_c)b,e\}}=r_ar_cW_{\{b,(r_cr_a)e\}}r_cr_a$. Hence
\begin{equation}\label{ineq}
(b+\sqrt2 c+\sqrt 2 a)\cdot e=(r_ar_c)b\cdot e>-1.
\end{equation}
Assume, for a contradiction, that $m_{ce}\ne 2$. Then $c\cdot e\le-1/2<1/2\sqrt 2$,
and so $(b+\sqrt 2 a)\cdot e>-1/2$, giving a contradiction if either $m_{be}\ne 2$ or
$m_{ae}\ne 2$. So $b\cdot e=a\cdot e=0$, and the inequality~\ref{ineq} above gives
$c\cdot e>1/\sqrt 2$. So $m_{ce}=3$. But now $W_{\{a,c,e\}}$ is of type $A_3$, hence
finite, and hence contained in a maximal finite subgroup that also
contains~$\FC(r_a)=W_J$. Since $b\in J$ it follows that $\{a,c,e,b\}$ is spherical,
which is false since it is of type $\widetilde B_3$. So $m_{ce}=2$. and it
remains to show that $m_{ae}=m_{be}=2$.

Since $c\cdot e=0$ we deduce from \ref{ineq} that $(b+\sqrt 2 a)\cdot e>-1$,
and in particular it follows that $m_{ae}$ is 2~or~3. In either case
$\{e,a,c\}$ is spherical (of type~$A_3$ or $A_1+A_2$), and so $\{e,a,c,b\}$ is
also spherical (since $r_b\in\FC(r_a)$). If either $m_{ae}\ne 2$ or $m_{be}\ne 2$
then applying Proposition~\ref{finclos2} with $L''=J\cup\{e,c\}$ in place of~$L$
yields a contradiction, since if $L_0''$ is
the the component of $L''$ containing~$a$ then $\{a,b\}\subseteq L_0''\cap J\ne L_0''$
(since $c\notin J$), but $L_0''$ is not of type $C_3$~or~$D_n$ since it
contains $\{a,c,b,e\}$. So $m_{ae}=m_{be}=2$, as required.
\end{proof}

\begin{lem}\label{alt1b}
Let $J'=J\setminus\{a\}$ and let $d\in\Odd(a)$. Then $m_{cd}\ne 2$.
If $m_{db'}\ne 2$ for some $b'\in J'$ then $m_{db'}=4$, and there
is a unique $a'$ adjacent to $d$ in $\Odd(a)$; moreover, $\{a',d,b'\}$
is of type~$C_3$, and $\FC(r_{a'})=W_{J'\cup\{a'\}}$. On the other
hand, if $m_{db'}=2$ for all $b'\in J'$ then $\FC(r_d)=W_{J'\cup\{d\}}$.
\end{lem}

\begin{proof}
We use induction on the distance from $d$ to $a$ in $\Odd(a)$. Observe that
if $d=a$ then $m_{db'}=2$ for all $b'\in J'$, since $\{a\}$ is a component
of~$J$, and we have $\FC(r_d)=W_J=W_{J'\cup\{d\}}$ and $m_{cd}=3\ne 2$, as
required.

Suppose now that $d\ne a$, and let $a=d_1,\,d_2,\,\ldots,\,d_k=d$ be a
minimal length path from $a$~to~$d$ in $\Odd(a)$. If $2\le i\le k-1$
then $d_i$ does not have valency~1 in $\Odd(a)$, and so $m_{d_ib'}=2$
for all $b'\in J'$, by the inductive hypothesis. The same is true for $i=1$,
since $\{a\}$ is a component of~$J$.

We prove first that $m_{cd}\ne 2$. Assuming, for a contradiction, that
$m_{cd}=2$, then clearly $d\notin\{a,c,b\}$, and Lemma~\ref{alt1a} tells
us that $m_{ad}=2$ and $m_{bd}=2$.
It follows that $m_{bf}=2$ for all $f$ in the set $M=\{d_1,d_2,\ldots,d_k\}$,
since $\{d_i\}$ is a component of $J'\cup\{d_i\}$ when $1\le i\le k-1$.
So $wb=b$ for all $w\in W_M$. Furthermore, since $d$~and~$a$ lie in the
same connected component of $\oddgraph(\Pi)$, we can choose $w\in W_M$
such that $wd=a$. Now since $wc\cdot a=c\cdot d=2$ we see that the
reflection $r_{wc}$ centralizes $r_a$, and hence normalizes $\FC(r_a)=W_J$.
By Lemma~\ref{subrootsys} it follows that either $wc\in\Phi_J$ or
$wc\cdot e=0$ for all $e\in J$. But $wc\cdot b=c\cdot w^{-1}b=c\cdot b\ne 0$;
so we must have $wc\in\Phi_J$, and hence
$c\in w^{-1}\Phi_J\subseteq\Phi_{J\cup M}$. So $c\in M$, contradicting $m_{bf}=2$
for all $f\in M$. So $m_{cd}\ne 2$.

Write $a'=d_{k-1}$ and $\widetilde J=J'\cup\{a'\}$. Note that $\{a'\}$ is a
component of $\widetilde J$, and $\FC(r_{a'})=W_{\widetilde J}$ 
(by the inductive
hypothesis). Now since $\{d,a'\}$ is spherical,
$\widetilde L=\widetilde J\cup\{d\}$ is spherical also. Let $\widetilde L_0$ be
the component of $\widetilde L$ containing~$a'$.

Consider first the case that $m_{db'}=2$ for all $b'\in J'$. Since also
$m_{a'b'}=2$ for all $b'\in J'$, it follows that $r_d$ and $r_{a'}$ both fix
all elements of~$J'$. Since $v=v[d,\{a'\}]\in W_{\{a',d\}}$ satisfies $va'=d$,
we conclude that
$$
\FC(r_d)=v\FC(r_{a'})v^{-1}=vW_{J'\cup\{a'\}}v^{-1}=W_{vJ'\cup\{va'\}}
=W_{J'\cup\{d\}},
$$
as required.

Now suppose that $m_{db'}\ne 2$ for some $b'\in J'$. Then
$b'\in \widetilde J\cap\widetilde L_0$, and so
$\{a'\}\ne\widetilde J\cap\widetilde L_0\ne\widetilde L_0$. Applying
Proposition~\ref{finclos2}, we see that the situation of alternative~(1)
must hold: alternative~(2) is ruled out since $a'$ is a component of
$\widetilde J$, and alternative~(3) is ruled out since $J'\cap\Odd(a)=\emptyset$.
Hence $\widetilde L_0=\{a',d,b'\}$ is of type~$C_3$, with $m_{db'}=4$ and
$m_{da'}=3$. Furthermore, Lemma~\ref{alt1a} tells us that
$m_{de}\in\{2,\infty\}$ for all~$e\in\Pi\setminus\{a',d,b'\}$; so $a'$
is the unique neighbour of $d$ in $\Odd(a)$, as required.
\end{proof}

\begin{lem}\label{alt1c}
Let $e\in\EOdd(a)\setminus\Odd(a)$ with $e\ne b$. Then $m_{be}=2$.
\end{lem}

\begin{proof}
If $e\in J$ then it is clear that $m_{be}=2$, since $b$ is a component
of~$J$. So we may assume that $e\notin J$.

As above, write $J'=J\setminus\{a\}$. Since $e\in\EOdd(a)$ there
exists a $d\in\Odd(a)$ with $m_{de}$ even. If $d$ is adjacent to some
$b'\in J'$ then, by Lemma~\ref{alt1b}, there is a unique $a'\in\Odd(a)$
adjacent to~$d$; furthermore, $\{a',d,b'\}$ is of type~$C_3$, and
$\FC(r_{a'})=W_{J'\cup\{a'\}}$. By Lemma~\ref{alt1a}, since
$m_{de}\ne\infty$ it follows that $m_{a'e}=m_{de}=2$. On the other
hand, if $d$ is not adjacent to any element of~$J'$ then Lemma~\ref{alt1b}
tells us that $\FC(r_d)=W_{J'\cup\{d\}}$. So in either case there is
an $a'\in\Odd(a)$ with $m_{a'e}$ even and $\FC(r_{a'})=W_{J'\cup\{a'\}}$.

Choose such an $a'$. Since $m_{a'e}$ is even, $v[e,\{a'\}]a'=a'$; moreover
$v[e,\{a'\}]$ is the reflection along some root $f=\lambda e+\mu a'$. Note
that $f\cdot a'=0$, and hence $\lambda\ne 0$. Since $e\notin J$ it follows
that $f\notin\Phi_J$. But $r_f$ centralizes $r_{a'}$, and hence normalizes
$\FC(r_{a'})=W_{J'\cup\{a'\}}$. By Lemma~\ref{subrootsys} it follows
that $f\cdot b=0$. But also $a'\cdot b=0$, since $\{a'\}$ and $\{b\}$
are distinct components of~$J'\cup\{a'\}$; so it follows that $e\cdot b=0$.
Thus $m_{be}=2$, as required.
\end{proof}

Lemmas \ref{alt1a}, \ref{alt1b} and \ref{alt1c} combine to yield the
following result.

\begin{prop}\label{alt1d}
Suppose that $a\in J\subseteq L\subseteq\Pi$, with $L$ spherical and
$\FC(r_a)=W_J$, and let $L_0$ be the component of~$L$ containing~$a$.
Suppose that $L_0=\{a,c,b\}$ is of type~$C_3$, with $m_{ac}=3$ and
$m_{bc}=4$, and $J\cap L_0=\{a,b\}$. Then $b$ is a $C_3$-neighbour
of~$\Odd(a)$. Furthermore, $J\cap\Odd(a)=\{a\}$, and if $a'\in\Odd(a)$
is not adjacent to any $C_3$-neighbour of $\Odd(a)$ then
$\FC(r_{a'})=W_{J'\cup\{a'\}}$, where $J'=J\setminus\{a\}$.
\end{prop}

\begin{proof}
If $m_{bd}\ne 2$ for some $d\in\Odd(a)$, then $m_{db}=4$, by
Lemma~\ref{alt1b}. There is at least one $d\in\Odd(a)$ such that $m_{bd}=4$,
namely $d=c$. Lemma~\ref{alt1b} tells us that for each $d\in\Odd(a)$ with
$m_{bd}=4$ there is an $a'\in\Odd(a)$ such that $\{a',d,b'\}$ is a system
of type~$C_3$. Moreover, by Lemma~\ref{alt1a}, if
$e\in\Pi\setminus(\Odd(a)\cup\{b\})$ then either $m_{de}=\infty$
or $m_{ae}=m_{be}=m_{de}=2$, while if $e\in\Odd(a)\setminus\{a,c\}$
then $m_{de}=\infty$, since $m_{de}\ne 2$ by Lemma~\ref{alt1b}. And
if $e\in\EOdd(a)\setminus(\Odd(a)\cup\{b\})$ then $m_{be}=2$, by
Lemma~\ref{alt1c}. So $b$ satisfies all the requirements of a
$C_3$-neighbour of $M=\Odd(a)$, as specified in
Definition~\ref{c3neighbour}.

It now follows from Lemma~\ref{alt1b} that if $a'\in\Odd(a)$ is adjacent
to some $b'\in J'$ then $b'$ is a $C_3$-neighbour of~$\Odd(a)$, and if
$a'$ is not adjacent to any such $b'$ then $\FC(r_{a'})=W_{J'\cup\{a'\}}$.
Finally, $J\cap\Odd(a)=\{a\}$, by Lemma~\ref{alt1a}.
\end{proof}

Let $a,\,a'\in\Pi$, and suppose that $w\in W$ has the property that
$wa=a'$. By Proposition~\ref{normalizers} there exist $a_i\in\Odd(a)$ and
$c_i\in\Pi$ such that
\begin{itemize}
\item[(i)]$a_1=a$ and $a_{k+1}=a'$,
\item[(ii)]$m_{c_ia_i}\ne\infty$ and $v[c_i,\{a_i\}]a_i=a_{i+1}$, for
all $i\in\{1,2,\ldots,k\}$,
\item[(iii)]$w=v[c_k,\{a_k\}]\cdots v[c_2,\{a_2\}]v[c_1,\{a_1\}]$.
\end{itemize}
Now let $b$ be a $C_3$-neighbour of~$\Odd(a)$. For each $c\in\Odd(a)$
that is adjacent to~$b$, define $X(c)=b+\sqrt 2c+\sqrt2\tilde a$, where
$\tilde a$ is the unique neighbour of $c$ in~$\Odd(a)$, and for each
$c\in\Odd(a)$ that is not adjacent to~$b$, define $X(c)=b$. We show that
$v[c_i,\{a_i\}]X(a_i)=X(a_{i+1})$ for all $i\in\{1,2,\ldots k\}$.

Suppose first that neither $c_i$ nor $a_i$ is adjacent to~$b$. Then
$X(a_i)=b$, and since $a_{i+1}\in\{a_i,c_i\}$ we have that $X(a_{i+1})=b$
also. Since $r_{a_i}$ and $r_{c_i}$ both fix $b$, and
$v[c_i,\{a_i\}]\in W_{\{a_i,c_i\}}$, it follows that
\[
v[c_i,\{a_i\}]X(a_i)=v[c_i,\{a_i\}]b=b=X(a_{i+1}),
\]
as required.

Next, suppose that $c_i$ is adjacent to~$b$, but $a_i$ is not adjacent to~$b$.
Since $m_{c_ia_i}\ne\infty$ and $a\in\Odd(a)$ it follows that $c_i\in\EOdd(a)$.
Since $b$ is a $C_3$-neighbour of $\Odd(a)$, it is not adjacent to any
element of $\EOdd(a)\setminus\Odd(a)$; so $c_i\in\Odd(a)$, and, moreover,
$m_{c_ie}=\infty$ for all $e\in\Odd(a)\setminus\{c_i\}$ apart from the
unique neighbour of $c_i$ in $\Odd(a)$. So $a_i$ is this unique neighbour,
$m_{c_ia_i}=3$, and $a_{i+1}=v[c_i,\{a_i\}]a_i=c_i$. Moreover,
$m_{c_ib}=4$ and $m_{a_ib}=2$. So
\[
v[c_i,\{a_i\}]X(a_i)=r_{a_i}r_{c_i}b=b+\sqrt2c_i+\sqrt2a_i=X(c_i)=X(a_{i+1})
\]
as required.

Now suppose that $a_i$ is adjacent to~$b$, and let $\tilde a$ be the unique
neighbour of $a_i$ in~$\Odd(a)$. Since $m_{a_ie}=\infty$ for all
$e\in\Odd(a)\setminus\{a_i,\tilde a\}$, if $c_i\in\Odd(a)$ then $c_i=\tilde a$.
In this case we see that
\[
v[c_i,\{a_i\}]X(a_i)=r_{a_i}r_{c_i}(b+\sqrt2c_i+\sqrt2a_i)=b=X(c_i)=X(a_{i+1}),
\]
since $a_{i+1}=v[c_i,\{a_i\}]a_i=c_i$. If $c_i=b$ then
$v[c_i,\{a_i\}]=r_br_{a_i}r_b$, which fixes both $a_i$ and
$X(a_i)=b+\sqrt 2a_i+\sqrt 2\tilde a$. So $v[c_i,\{a_i\}]X(a_i)=X(a_{i+1})$
in this case too. Finally, suppose that $c_i\notin\Odd(a)\cup\{b\}$.
Since $m_{c_ia_i}\ne\infty$ we must have $m_{c_i\tilde a}=m_{c_ia_i}=m_{c_ib}=2$,
(by the definition of a $C_3$-neighbour). So
\[
v[c_i,\{a_i\}]X(a_i)=r_{c_i}(b+\sqrt2c_i+\sqrt2\tilde a)
=b+\sqrt2c_i+\sqrt2\tilde a=X(a_{i+1})
\]
since $a_{i+1}=r_{c_i}a_i=a_i$.

We have now now covered all cases, and shown that
$v[c_i,\{a_i\}]X(a_i)=X(a_{i+1})$ for all $i\in\{1,2,\ldots k\}$. By a trivial
induction it follows that $X(a_{k+1})=wX(a_1)$.

Thus we have established the following result.

\begin{lem}\label{c3lem}
Let $a\in\Pi$ and $w\in W$ such that $wa\in\Pi$. Suppose that $b$ is a
$C_3$-neighbour of $\Odd(a)$ that is not adjacent to~$a$. Then
\[
wb=\begin{cases}
b&\text{if $wa$ is not adjacent to~$b$}\\
b+\sqrt2 wa+\sqrt2\tilde a&\text{if $wa$ is adjacent to~$b$}
\end{cases}
\]
where $\tilde a$ is adjacent to $wa$ in~$\Odd(a)$.
\end{lem}

We are now able to give a detailed description of the components of $J$
whenever $W_J$ is the finite continuation of a simple reflection.

\begin{prop}\label{components}
Suppose that $a\in J\subseteq\Pi$ with $W_J=\FC(r_a)$, and suppose that $K$ is
a component of~$J$. Then one of the following alternatives holds.
\begin{itemize}
\item[(a)] $K=\{a\}=J\cap\Odd(a)$.
\item[(b)] $K=\{a,b\}$ is of type~$C_2$, and $J\cap\Odd(a)=\{a\}$.
\item[(c)] $K=\{a\}$ or $K=\{b\}$, where $\{a,b\}=J\cap\Odd(a)$ is of type $A_1+A_1$.
\item[(d)] $K=\{b\}\nsubseteq\Odd(a)$, and $b$ is a $C_3$-neighbour of $\Odd(a)$.
\item[(e)] $\Odd(a)\subseteq K$, and $K$ is a component of~$\EOdd(a)$.
\item[(f)] $K\cap\Odd(a)=\emptyset$, and $K$ is a component of~$\EOdd(a)$.
\end{itemize}
\end{prop}

\begin{proof}
We consider first the case that $K\cap\Odd(a)\ne\emptyset$, and start by
supposing that there exists a spherical $L\subseteq\Pi$ with $J\subseteq L$
and $K$ not a component of~$L$.

Choose such an~$L$, and let $L_0$ be the component of $L$ containing~$a$.
By Proposition~\ref{finclos2}, since $K$ is not a component of $L$ we must
have $K\subseteq L_0$. So either $K=\{a\}$, in which case (a) above holds,
or else $\{a\}\subsetneqq \{a\}\cup K\subseteq J\cap L_0$. Furthermore,
$J\cap L_0\ne L_0$, since $K\ne L_0$. So if (a) does not hold then
$\{a\}\ne J\cap L_0\ne L_0$, and so one of the alternatives (1), (2) or~(3)
of Proposition~\ref{finclos2} must hold.

Suppose that alternative (2) holds, so that $K=\{a,b\}=J\cap L_0$
for some~$b$, and $\{a,b\}$ is of type~$C_2$. By Proposition~\ref{alt2} we
see that each $c\in\Odd(a)\setminus\{a\}$ lies in a type $C$ spherical
subset $L'$ of $\Pi$ containing~$\{a,b\}$. Since $J\cap L'=\{a,b\}$ (by
Proposition~\ref{finclos2}) it follows that $c\notin J$. So
$J\cap\Odd(a)=\{a\}$, and (b) above is satisfied.

Suppose that alternative (3) of Proposition~\ref{finclos2} holds, so
that $J\cap L_0=\{a,b\}$ is of type $A_1+A_1$, and $b\in\Odd(a)$.
Proposition~\ref{alt3} immmediately yields that $J\cap\Odd(a)=\{a,b\}$,
and so (c) above is satisfied.

Suppose that alternative~(1) of Proposition~\ref{finclos2} holds, so that
$L_0=\{a,c,b\}$ with $m_{ac}=3$ and $m_{cb}=4$, and $J\cap L_0=\{a,b\}$.
By Lemma~\ref{alt1a} we know that $b\notin\Odd(a)$, and since we have
assumed that $K\cap\Odd(a)\ne\emptyset$, it follows that
$K=\{a\}=J\cap\Odd(a)$. So (a) holds.

We have now dealt with all cases that arise if there is a spherical
$L\subseteq\Pi$ with $J\subseteq L$ and $K$ not a component of~$L$.
So assume that $K$ is a component of every spherical $L$ containing~$J$.
We show that in this case $\Odd(a)\subseteq K$, and $K$ is a component
of~$\EOdd(a)$; that is, (e) above holds.

To show that $\Odd(a)\subseteq K$ it is clearly sufficient to show
that if $a'\in K\cap \Odd(a)$ and $b$ is adjacent to $a'$ in $\Odd(a)$
then $b\in K$. Note that since $a'\in\Odd(a)$ there exists $w\in W$
with $a'=wa$, and since Proposition~\ref{finclos} yields that
$w\in N_W(W_J)$ it follows that $\FC(r_{a'})=W_J$. Now the assumption
that $b$ and $a'$ are adjacent in $\Odd(a)$ implies that $\{a',b\}$ is
spherical, and therefore $J\cup\{b\}$ is spherical. But $K$ is a
component of every spherical subset of $\Pi$ containing~$J$; so it is a
component of $J\cup\{b\}$. But $a'\in K$ and $b$ is adjacent to $a'$; so
$b\in K$, as required.

Since $K\subseteq J\subseteq\EOdd(a)$ and $K$ is connected, saying that
$K$ is a component of $\EOdd(a)$ is equivalent to saying
that $m_{bc}=2$ whenever $b\in K$ and $c\in\EOdd(a)\setminus K$.
So suppose that $c\in \EOdd(a)\setminus K$. Then there exists an
$a'\in\Odd(a)$ such that $m_{a'c}$ is even. Thus $\{a',c\}$ is
spherical, and as above it follows that $J\cup\{c\}$ is spherical.
So $K$ must be a component of $J\cup\{c\}$, and since $c\notin K$
it follows that $m_{bc}=2$ for all $b\in K$, as required.

%
%
%

It remains to consider the case that $K\cap\Odd(a)=\emptyset$; we must
show that either (f) or (d) holds. We start
by supposing that there exists a spherical $L\subseteq\Pi$ and a $w\in W$
with $wJ\subseteq L$ and $wK$ not a component of~$L$.

Choose such $L$~and~$w$, and let $L_0$ be the component of $L$
containing~$wa$. By Proposition~\ref{finclos2}, since $wK$ is not a
component of $L$ we must have $wK\subseteq L_0$. Now $wJ\cap L_0\ne L_0$
since $wK\ne L_0$, and $\{wa\}\ne wJ\cap L_0$ since $wa\notin wK$.
So one of the alternatives (1), (2) or~(3) of Proposition~\ref{finclos2}
must hold. Alternative~(3) can be ruled out, since in that case
$wJ\cap L_0\subseteq\Odd(wa)$, which is impossible since
$K\cap\Odd(a)=\emptyset$. If alternative~(2) holds then
$wK=wJ\cap L_0$ contains $wa$ and is of type~$C_2$, whence $K$
contains $a$ and is of type~$C_2$, and (b) is satisfied. If
alternative~(1) holds then since $wK\ne\{wa\}$ it follows from
Proposition~\ref{alt1c} that $wK=\{b\}$, with $b$ a $C_3$-neighbour of
$\Odd(a)$. Since $wa$ is not adjacent to $b$, it follows from
Lemma~\ref{c3lem} that $w^{-1}b=b$, unless $a$ is adjacent to~$b$,
in which case $w^{-1}b=b+\sqrt 2a+\sqrt 2\tilde a$ for some $\tilde a$
in $\Odd(a)$. But this latter case cannot occur, since
$w^{-1}b\in K\subseteq\Pi$. So $K=wK=\{b\}$, with $b$ a $C_3$-neighbour
of $\Odd(a)$, and (d) holds.

Finally, suppose that $wK$ is a component of every spherical
$L\subseteq\Pi$ such that $wJ\subseteq L$ for some $w\in W$.
For each $c\in\EOdd(a)\setminus K$ there is then a sequence
$a=a_0,a_1,\ldots,a_k=c$ in $\Pi$ such that
$m_{a_{i-1}a_i}$ finite for all $i\in\{1,2,\ldots,k\}$ and odd for all
$i\in\{1,2,\ldots,k-1\}$. We shall show that, for every such sequence,
$m_{ba_i}=2$ for all $b\in K$ and $i\in\{0,1,\ldots,k\}$; in particular,
this will show that $m_{bc}=2$ whenever $b\in K$ and $c\in\EOdd(a)\setminus K$,
enabling us to conclude that $K$ is a component of~$\EOdd(a)$.

The case $k=0$ is clear, since $a\in J\setminus K$ and $K$ is a component of~$J$.
Proceeding by induction, we may assume that $k>0$ and $m_{ba_i}=2$ for all
$i\in\{1,2,\ldots,k-1\}$ and all $b\in K$. We see that the element
$u=v[a_{k-1},\{a_{k-2}\}]v[a_{k-2},\{a_{k-3}\}]\cdots v[a_1,\{a_0\}]$ centralizes
$W_K$ and has the property that $ua=a_{k-1}$, since the labels in the path from
$a$ to $a_{k-1}$ are all odd. The group $u^{-1}W_{\{a_{k-1},a_k\}}u$ is finite
and contains $u^{-1}r_{a_{k-1}}u=r_a$, and so there is a maximal finite subgroup
$G$ of $W$ containing this group and also containing~$W_J$.

Note that $W_J\cup\{u^{-1}r_cu\}\subseteq G=w^{-1}W_Lw$, for some $w\in W$ and
spherical $L\subseteq\Pi$, the element $u$ being in the centralizer of~$W_K$. We
may choose $w$ to be the minimal length element of $W_Lw=W_LwW_J$, and it
then follows from Lemma~\ref{kilmoyerlem} that $wJ\subseteq L$. Hence $wK$ is a
component of~$L$. Furthermore, since $wu^{-1}r_cuw^{-1}\in W_L$ we see that
the root $wu^{-1}c$ is in $\Phi_L$ and not in $\Phi_{wK}=wu^{-1}\Phi_K$ (since
$c\notin\Phi_K$). So $wu^{-1}c\cdot wu^{-1}b=0$ for all $b\in K$.
So $c\cdot b=0$, or (equivalently) $m_{bc}=2$ for all $b\in K$, as required.
\end{proof}

To complement the results we have obtained so far, our next task is to find
conditions that ensure that a visible subgroup~$W_K$ is contained in $\FC(r_a)$.

\begin{lem}\label{lem1a} Let $a\in\Pi$ and $K$ a component
of $\EOdd(a)$ such that $W_K$ is finite. Then $W_K\subseteq \FC(r_a)$.
\end{lem}

\begin{proof}
Let $F$ be a maximal finite subgroup of $W$ with
$r_a\in F$, and choose $w\in W$ such that $wFw^{-1}=W_L$ for some
$L\subseteq\Pi$. We may replace $w$ by the minimal length element
in the double coset $W_LwW_{\{a\}}$, since this does not affect the
condition $wFw^{-1}=W_L$. So we have that $w^{-1}L\subseteq\Phi^+$,
and, moreover, $r_a\in w^{-1}Lw\cap W_{\{a\}}=W_{w^{-1}L\cap\{a\}}$
by Lemma~\ref{kilmoyerlem}. So $wa\in L\subseteq\Pi$,
and by Lemma~\ref{normalizers} we see that $w$ is a product of factors of
the form $v[d,\{c\}]$, with $c,\,d\in\EOdd(a)$. Since $K$ is a
component of $\EOdd(a)$ it follows that each $v[d,\{c\}]$ normalizes~$W_K$,
and therefore $w$ normalizes~$W_K$. Moreover, since $wa\in L$ and $L$ is
spherical, it follows that $L\subseteq\EOdd(wa)=\EOdd(a)$. So $W_L$
normalizes~$W_K$. But $W_K$ is finite, by hypothesis, and $W_L$ is a
maximal finite subgroup of~$W$. So $W_K\subseteq W_L$, and
$W_K=w^{-1}W_Kw\subseteq w^{-1}W_Lw=F$. Thus $W_K$ is contained
in all maximal finite subgroups of $W$ containing~$r_a$, as required.
\end{proof}

\begin{lem}\label{morec3}
Let $a\in\Pi$ and let $b$ be a $C_3$-neighbour of $\Odd(a)$. If
$a$~and~$b$ are not adjacent in $\Pi$ then $r_b\in\FC(r_a)$.
\end{lem}

\begin{proof}
Let $F$ be a maximal finite subgroup of $W$ with $r_a\in F$. As in the
proof of Lemma~\ref{lem1a} there exist a $w\in W$ and a
maximal spherical $L\subseteq\Pi$ with $wa=a'\in L$ and $F=w^{-1}W_Lw$.
Since $L$ is spherical, $L\subseteq\EOdd(a)$.

Suppose first that $a'$ is not adjacent to~$b$. Then $m_{ca'}=\infty$
for every $c\in\Odd(a)$ that is adjacent to~$b$, and since $a'\in L$ it follows
that no such $c$ is in~$L$. Thus $m_{be}=2$ for all $e\in L\cap\Odd(a)$.
But since also $m_{be}=2$ for all $e\in\EOdd(a)\setminus(\Odd(a)\cup\{b\})$, it
follows that $m_{be}=2$ for all $e\in L\setminus\{b\}$. Thus $\{b\}$ is
a component of $L\cup\{b\}$, and since $L$ is spherical it follows that
$L\cup\{b\}$ is spherical. Maximality of $L$ tells us that $b\in L$.
Moreover, Lemma~\ref{c3lem} gives $wb=b$, and so
$r_b=w^{-1}r_bw\in w^{-1}W_Lw=F$.

On the other hand, suppose that $a'$ is adjacent to~$b$. In this case
Lemma~\ref{c3lem} gives $wb=b+\sqrt 2a'+\sqrt 2\tilde a$, where $\tilde a$
is the unique neighbour of $a'$ in~$\Odd(a)$. Furthermore, since
$m_{a'e}\in\{2.\infty\}$ for all $e\in\Pi\{\tilde a,a',b\}$, we see that
$m_{a'e}=2$ for all $e\in L\setminus\{\tilde a,a',b\}$ (since $L$ is
spherical). But the definition of a $C_3$-vertex also requires that
$m_{\tilde ae}=m_{be}=2$ whenever $m_{a'e}=2$; so it follows that
$\{\tilde a,a',b\}$ is a component of $L\cup\{\tilde a,a',b\}$, which
is therefore spherical since $L$ and $\{\tilde a,a',b\}$ are both
spherical. Maximality of $L$ tells us that $\{\tilde a,a',b\}\subseteq L$;
so $wb=b+\sqrt 2a'+\sqrt 2\tilde a\in\Phi_L$, and
$r_b=w^{-1}r_{wb}w\in w^{-1}W_Lw=F$.

So $r_b\in F$ in all cases, and so $r_b$ is contained
in all maximal finite subgroups of $W$ containing~$r_a$, as required.
\end{proof}

We now prove the converse to Proposition~\ref{alt2}.

\begin{prop}\label{alt2+}
Let $a\in\Pi$ and $b\in\Pi\setminus\Odd(a)$, and suppose that $(a,b)$ is a focus
of $\Odd(a)$ in $\Pi$. Then $\FC(r_a)=W_J$, where $J$ is the union of
$\{a,b\}$ and the spherical components of $\EOdd(a)$. Moreover, $\FC(r_{a'})$
is not visible for any $a'\in\Odd(a)\setminus\{a\}$.
\end{prop}

\begin{proof}
For each $c\in\Odd(a)$ let $X(c)=b+\sqrt 2\sum_{i=1}^m c_i$
and $Y(c)=b+\sqrt 2\sum_{i=1}^{m-1} c_i$, where $c_1=a,\,c_2,\,\ldots,\,c_m=c$
is the unique path from $a$~to~$c$ in $\Odd(a)$, noting that $X(c)$ and $Y(c)$
are roots in $\Phi_{C[b..c]}$. We remark, for later use, that $X(c)$ and $Y(c)$
are fixed by the reflections $r_b,\,r_{c_1},\,\ldots,\,r_{c_{m-2}}$.

Let $F=w^{-1}W_Lw$ be a maximal finite subgroup of $W$ containing~$r_a$, with
$L\subseteq\Pi$ and $w$ of minimal length in~$W_LwW_{a}$. Then $wa=a'\in L$,
by Lemma~\ref{kilmoyerlem}. Put $L_0=L\cap\Odd(a)$.

Choose $c\in L_0$ with $C[b..c]$ of maximal cardinality. If $d\in L_0$ then
$m_{cd}\ne\infty$ (since $L_0$ is spherical), whence $d\in C[b..c]$ by
condition~(3) of Definition~\ref{focus}. So $L_0\subseteq C[b..c]$.
Now if $e\in L\setminus L_0$ is arbitrary then $e\notin\Odd(a)$ (since
$e\notin L\cap\Odd(a)$) and $m_{ce}<\infty$ (since $c,\,e\in L$ and $L$ is
spherical). By condition~(4) of Definition~\ref{focus}
it follows that $m_{de}=2$ for all
$d\in C[b..c]$. Since this holds for all $e\in L\setminus L_0$, and
$C[b..c]$ and $L\setminus L_0$ are both spherical, it follows that
$C[b..c]\cup(L\setminus L_0)$ is spherical. But this set contains $L$ (since
$L_0\subseteq C[b..c]$) and since $L$ is a maximal spherical subset of~$\Pi$
we conclude that $L=C[b..c]\cup(L\setminus L_0)$.

By Proposition~\ref{normalizers} and Lemma~\ref{class} there exist simple roots
$e_1,\,e_2,\,\ldots,\,e_k$ and $d_1=a,\,d_2,\,\ldots,\,d_{k+1}=a'\in\Odd(a)$
with $w=v[e_k,\{d_k\}]\cdots v[e_2,\{d_2\}]v[e_1,\{d_1\}]$
and $v[e_i,\{d_i\}]d_i=d_{i+1}$ for all $i\in\{1,2,\ldots,k\}$. Moreover,
$m_{e_id_i}<\infty$ for all~$i$. Let $w_0=1$ and $w_i=v[e_i,\{d_i\}]w_{i-1}$;
we will show that
\[
\{w_ib,-w_ib,w_i(b+\sqrt 2a),-w_i(b+\sqrt 2a)\}
=\{X(d_{i+1}),-X(d_{i+1}),Y(d_{i+1}),-Y(d_{i+1})\}
\]
for all $i\in\{0,1,\ldots,k\}$. The case $i=0$ is trivial.

Proceeding inductively, suppose that $i>1$ and
\[
\{\pm w_{i-1}b,\pm w_{i-1}(b+\sqrt2a)\}=\{\pm X(d_i),\pm Y(d_i)\}.
\]
It will be sufficient to show that $v[e_i,\{d_i\}]X(d_i)$ and
$v[e_i,\{d_i\}]Y(d_i)$ both lie in the set $\{\pm X(d_{i+1}),\pm Y(d_{i+1})\}$.

Suppose first that $d_i=a$. Then $X(d_i)=b+\sqrt2 a$ and $Y(d_i)=b$.
If $e_i\notin\Odd(a)$ then $m_{e_id_i}$ is even, and $d_{i+1}=d_i=a$.
Furthermore, by condition (4) of Definition~\ref{focus}we
have either $\{e_i,d_i\}=\{b,a\}$ or $m_{e_ib}=m_{e_ia}=2$. In the
former case  $v[e_i,\{d_i\}]=v[b,\{a\}]=r_br_ar_b$, giving
$v[e_i,\{d_i\}]b=-b-\sqrt 2a=-X(a)$
and $v[e_i,\{d_i\}](b+\sqrt2a)=-b=-Y(a)$; in the latter case
$v[e_i,\{d_i\}]=v[e_i\{a\}]=r_{e_i}$, giving $v[e_i,\{d_i\}]b=b=Y(a)$ and
$v[e_i,\{d_i\}](b+\sqrt 2a)=b+\sqrt2a=X(a)$. If $e_i\in\Odd(a)$ then
$a\in C[b..e_i]$, and by condition~(2) of Definition~\ref{focus} we have
$m_{e_ib}=2$ and either $m_{e_ia}=2$ or $m_{e_ia}=3$. If $m_{e_ia}=2$ then
$d_{i+1}=d_i=a$, while if $m_{e_ia}=3$
then $d_{i+1}=e_i$. Furthermore, in former case we find that
$v[e_i,\{d_i\}]b=r_{e_i}b=b=Y(a)$ and $v[e_i,\{d_i\}](b+\sqrt2a)=b+\sqrt2a=X(a)$,
while in the latter case we find that
$v[e_i,\{d_i\}]b=r_ar_{e_1}b=b+\sqrt 2a=Y(e_i)$ and
$v[e_i,\{d_i\}](b+\sqrt2a+\sqrt2e_i)=X(e_i)$.

Now suppose that $d_i\ne a$. If $e_i\notin\Odd(a)\cup\{b\}$ then
$m_{e_id_i}=2$ and $d_{i+1}=d_i$. Moreover, $m_{e_id}=2$ for all
$d\in C[b..d_i]$, and so $v[e_i,\{d_i\}]=r_{e_i}$ fixes all the roots
in $\Phi_{C[b..d_i]}$, including $X(d_i)=X(d_{i+1})$ and $Y(d_i)=Y(d_{i+1})$.
If $e_i\in\Odd(a)\cup\{b\}$ and $\{e_i,d_i\}$ is not an edge of $\Odd(a)$
then we again have $d_{i+1}=d_i$ and $v[e_i,\{d_i\}]=r_{e_i}$. By
condition~(3) of Definition~\ref{focus} we either have
$d_i\in C[b..e_i]$ or $e_i\in C[b..d_i]$. In the former case we have
$m_{e_id}=2$ for all~$d\in C[b..d_i]$, and as above we see that $r_{e_i}$
fixes $X(d_i)$ and $Y(d_i)$. In the latter case the remark made at the
start of the proof implies that it is still true that $r_{e_i}$
fixes $X(d_i)$ and $Y(d_i)$. So we have shown that when $m_{e_id_i}=2$
it is true that $v[e_i,\{d_i\}]X(d_i)$ and $v[e_i,\{d_i\}]Y(d_i)$ both lie in
the set $\{\pm X(d_{i+1}),\pm Y(d_{i+1})\}$, and it remains only to consider the
case that $e_i$ and $d_i$ are adjacent in~$\Odd(a)$. Note that in this case
$d_{i+1}=e_i$.

Let $C[b..d_i]=\{b,c_1,\ldots,c_m\}$ with $c_1=a$ and $c_m=d_i$, and suppose
that $e_i=c_{m-1}$ is the vertex adjacent to~$d_i$ in $C[b..d_i]$. Then
\[
\displaylines{
v[e_i,\{d_i\}]X(d_i)=r_{c_m}r_{c_{m-1}}\Bigl(b+\sqrt2\sum_{j=1}^mc_j\Bigr)
=r_{c_m}\Bigl(b+\sqrt2\sum_{j=1}^mc_j\Bigr)\cr
=b+\sqrt2\sum_{j=1}^{m-1}c_j=X(e_i),\cr}
\]
and similarly
\[
\displaylines{
v[e_i,\{d_i\}]Y(d_i)=r_{c_m}r_{c_{m-1}}\Bigl(b+\sqrt2\sum_{j=1}^{m-1}c_j\Bigr)
=r_{c_m}\Bigl(b+\sqrt2\sum_{j=1}^{m-2}c_j\Bigr)\cr
=b+\sqrt2\sum_{j=1}^{m-2}c_j=Y(e_i).\cr}
\]
The alternative possibility is that $d_i$ is adjacent to $e_i$ in
$C[b..e_i]$. Exactly the same calculations show that
$v[e_i,\{d_i\}]X(d_i)=X(e_i)$ and $v[e_i,\{d_i\}]Y(d_i)=Y(e_i)$ in this case also.

The induction is now complete, and it follows in particular that
$wb=w_kb$ is one of $\pm X(a')$ or $\pm Y(a')$. Hence
\[
wb\in\Phi_{C[b..a']}\subseteq\Phi_{C[b..c]}\subseteq\Phi_L.
\]
Thus $wr_bw^{-1}\in W_L$, and so $r_b\in w^{-1}W_Lw=F$. Since $F$ was an arbitrary
maximal finite subgroup of~$W$ containing~$r_a$, this shows that $r_b\in\FC(r_a)$.

Let $\widetilde M$ be the component of $\EOdd(a)$ containing $\Odd(a)$, and suppose,
for a contradiction, that $\widetilde M$ is spherical. Clearly $b\in\widetilde M$,
since $m_{ba}=4$, but $\widetilde M=\Odd(a)\cup\{b\}$ is not permitted, in
view of condition~(5) of Definition~\ref{focus}. So
$\widetilde M\setminus(\Odd(a)\cup\{b\}\ne\emptyset$. But for
$e\in\widetilde M\setminus(\Odd(a)\cup\{b\})$ and $c\in\Odd(a)$ we have
$m_{ce}\ne\infty$, since $\widetilde M$ is spherical, and by condition~(4)
of~Definition~\ref{focus} it follows that $m_{be}=m_{ce}=2$ for all
$c\in\Odd(a)$. This contradicts the fact that $\widetilde M$ is connected.

Now suppose that $a'\in\Odd(a)$ is such that $\FC(r_{a'})=W_J$ for some
$J\subseteq\Pi$, and let $J_0$ be the component of $J$ containing~$a'$.
Since $J_0\ne\widetilde M$ it follows from Proposition~\ref{components}
that $J_0$ has rank at most~2. Now since there
exists $w\in W_{C[b..a']}$ such that $wa=-a'$ and $wb=X(a')$, and since
$r_b\in\FC(r_a)$, it follows that $r_{wb}=wr_bw^{-1}\in\FC(wr_aw^{-1})=\FC(r_{a'})$.
Thus $X(a')\in\Phi_J$, and so $C[b..a']\subseteq J$. Since $J_0$ has rank
at most~2, this means that $a'=a$ and $J_0=\{a,b\}$.

It remains to prove that $J$ is the union of $J_0$ and the spherical components
of~$\EOdd(a)$. By Lemma~\ref{lem1a} we know that all these components are
contained in~$J$. But if $K$ is any other component of $J$ such that
$K\cap\Odd(a)=\emptyset$,
then by Proposition~\ref{components} we see that $K=\{b'\}$, with $b'$ a
$C_3$-neighbour of~$\Odd(a)$. Since $b$ is the only element of $\Pi$ such that
$m_{bc}\in\{2,4\}$ for all $c\in\Odd(a)$, we must have $b'=b$, contradicting
the fact that the component of $J$ containing~$b$ is $J_0=\{a,b\}$.
\end{proof}

Next, we have the converse to Proposition~\ref{alt3}.

\begin{prop}\label{alt3+}
Let $a\in\Pi$ and suppose that there exists a $b\in\Odd(a)$ such that
$\{a,b\}$ is a half-focus of $\Odd(a)$ in~$\Pi$. Suppose also that the
vertices $\Odd(a)$ do not comprise a spherical subset of~$\Pi$.
Then $\FC(r_a)=W_J$, where $J$ is the union of $\{a,b\}$ and the spherical
components of $\EOdd(a)$. Moreover, $\FC(r_{a'})$
is not visible for any $a'\in\Odd(a)\setminus\{a,b\}$.
\end{prop}

%

\begin{proof}
For each $c\in\Odd(a)\setminus\{a,b\}$, define
\[
X(c)=b+a+c+2\sum_{i=2}^{m-1}c_i
\]
where $c_1=a,\,c_2,\,\ldots,\,c_m=c$ is the unique path from $a$~to~$c$ in
$\Odd(a)\setminus\{b\}$. Then $X(c)$ is a root in $\Phi_{D[a,b..c]}$ and is
fixed by the reflections $r_b,\,r_{c_1},\,\ldots,\,r_{c_{m-2}}$ and~$r_{c_m}$.
Define also $X(a)=b$ and $X(b)=a$.

Let $F=w^{-1}W_Lw$ be a maximal finite subgroup of $W$ containing~$r_a$, with
$L\subseteq\Pi$ and $w$ of minimal length in~$W_LwW_{a}$. Then $wa=a'\in L$,
by Lemma~\ref{kilmoyerlem}. Put $L_0=L\cap\Odd(a)$.

Choose $c\in L_0$ with $D[a,b..c]$ of maximal cardinality. If $d\in L_0$ then
$m_{cd}\ne\infty$ (since $L_0$ is spherical), whence $d\in D[a,b..c]$
by condition~(4) of Definition~\ref{halffocus}. So
$L_0\subseteq D[a,b..c]$.
Now if $e\in L\setminus L_0$ is arbitrary then $e\notin\Odd(a)$ (since
$e\notin L\cap\Odd(a)$) and $m_{ce}<\infty$ (since $c,\,e\in L$ and $L$ is
spherical). By condition~(5) of Definition~\ref{halffocus} it follows that
$m_{de}=2$ for all
$d\in D[a,b..c]$. Since this holds for all $e\in L\setminus L_0$, and
$D[a,b..c]$ and $L\setminus L_0$ are both spherical, it follows that
$D[a,b..c]\cup(L\setminus L_0)$ is spherical. But this set contains $L$ (since
$L_0\subseteq D[a,b..c]$) and since $L$ is a maximal spherical subset of~$\Pi$
we conclude that $L=D[a,b..c]\cup(L\setminus L_0)$.

By Proposition~\ref{normalizers} and Lemma~\ref{class} there exist simple roots
$e_1,\,e_2,\,\ldots,\,e_k$ and $d_1=a,\,d_2,\,\ldots,\,d_{k+1}=a'\in\Odd(a)$
with $w=v[e_k,\{d_k\}]\cdots v[e_2,\{d_2\}]v[e_1,\{d_1\}]$
and $v[e_i,\{d_i\}]d_i=d_{i+1}$ for all $i\in\{1,2,\ldots,k\}$. Furthermore,
we have $m_{e_id_i}<\infty$ for all~$i$. Let $w_0=1$, and
$w_i=v[e_i,\{d_i\}]w_{i-1}$ for $i\ge 1$. We will show that
\[
\{w_ib,-w_ib\}=\{X(d_{i+1}),-X(d_{i+1})\}
\]
for all $i\in\{0,1,\ldots,k\}$.

The case $i=0$ is trivial.
Proceeding inductively, suppose that $i>1$ and $w_{i-1}b=\pm X(d_i)$.
It will be sufficient to show that $v[e_i,\{d_i\}]X(d_i)=\pm X(d_{i+1})$.

Suppose first that $d_i=a$, so that $X(d_i)=b$. If $e_i\ne b$ then
$m_{e_ib}=m_{e_ia}\in\{2,3\}$, since $m_{e_ia}=m_{e_id_i}\ne\infty$.
We also have $m_{e_ia}=2$ if $e_i=b$. In the case $m_{e_ia}=3$ we
have $v[e_i,\{d_i\}]=r_ar_{e_i}$, and $d_{i+1}=r_ar_{e_i}a=e_i$.
Furthermore,
\[
v[e_i,\{d_i\}]X(d_i)=r_ar_{e_i}b=a+b+e_i=X(e_i)=X(d_{i+1}),
\]
as required. In the case $m_{e_ia}=2$ we have $v[e_i,\{d_i\}]=r_{e_i}$,
giving $d_{i+1}=r_{e_i}a=a$, and
\[
v[e_i,\{d_i\}]X(d_i)=r_{e_i}b=\pm b=\pm X(d_{i+1}),
\]
since either $e_i=b$ or $m_{e_ib}=2$.

The case $d_i=b$ is the same as the case $d_i=a$ with $a$ and $b$
interchanged; so suppose that $d_i\notin\{a,b\}$.
If $e_i\notin\Odd(a)$ then $m_{e_id_i}=2$ and $d_{i+1}=d_i$. Moreover,
$m_{e_id}=2$ for all $d\in D[a,b..d_i]$, and so $v[e_i,\{d_i\}]=r_{e_i}$
fixes all the roots in $\Phi_{D[a,b..d_i]}$, including $X(d_i)=X(d_{i+1})$.
If $e_i\in\Odd(a)$ and $\{e_i,d_i\}$ is not an edge of $\Odd(a)$
then we again have $d_{i+1}=d_i$ and $v[e_i,\{d_i\}]=r_{e_i}$. By
condition~(4) of Definition~\ref{halffocus} we either have
$d_i\in D[a,b..e_i]$ or $e_i\in D[a,b..d_i]$.
In the former case we have $m_{e_id}=2$ for all~$d\in C[b..d_i]$, and as
above we see that $r_{e_i}$ fixes $X(d_i)$. In the latter case it is still
true that $r_{e_i}$ fixes $X(d_i)$, since the only simple reflection of
$D[a,b..d_i]$ that does not fix $X(d_i)$ is the one corresponding to the vertex
adjacent to~$d_i$. So we have shown that when $m_{e_id_i}=2$
it is true that $v[e_i,\{d_i\}]X(d_i)$ and $v[e_i,\{d_i\}]Y(d_i)$ both lie in
the set $\{\pm X(d_{i+1}),\pm Y(d_{i+1})\}$, and it remains to consider the
case that $e_i$ and $d_i$ are adjacent in~$\Odd(a)$. Note that in this
case~$d_{i+1}=e_i$.

Let $D[a,b..d_i]=\{b,c_1,\ldots,c_m\}$ with $c_1=a$ and $c_m=d_i$. Suppose first
that $m>2$, and suppose that $e_i=c_{m-1}$ is the vertex adjacent to~$d_i$ in
$D[a,b..d_i]$. Then
\begin{align*}
v[e_i,\{d_i\}]X(d_i)&=r_{c_m}r_{c_{m-1}}\Bigl(b+a+c_m+2\sum_{j=2}^{m-1}c_j\Bigr)\\
&=r_{c_m}\Bigl(b+a+c_m+c_{m-1}+2\sum_{j=1}^{m-2}c_j\Bigr)\\
&=b+a+c_{m-1}+2\sum_{j=1}^{m-2}c_j=X(e_i).
\end{align*}
If $m=2$ and $e_i=b$ then
\[
v[e_i,\{d_i\}]X(d_i)=r_br_{d_i}(a+b+d_i)=a=X(b)=X(e_i),
\]
and if $e_i=a$ then similarly
\[
v[e_i,\{d_i\}]X(d_i)=r_ar_{d_i}(a+b+d_i)=b=X(a)=X(e_i).
\]
The alternative possibility is that $d_i=c_{m-1}$ is the vertex adjacent to
$e_i=c_m$ in $D[a,b..e_i]=\{b,c_1,\ldots,c_m\}$. We calculate that
\begin{align*}
v[e_i,\{d_i\}]X(d_i)
&=r_{c_{m-1}}r_{c_m}\Bigl(b+a+c_{m-1}+2\sum_{j=2}^{m-2}c_j\Bigr)\\
&=r_{c_{m-1}}\Bigl(b+a+c_m+c_{m-1}+2\sum_{j=1}^{m-2}c_j\Bigr)\\
&=b+a+c_m+2\sum_{j=1}^{m-1}c_j=X(e_i),
\end{align*}
as required.

The induction is now complete, and it follows that $wb=w_kb=\pm X(a')$. Hence
\[
wb\in\Phi_{D[a,b..a']}\subseteq\Phi_{D[a,b..c]}\subseteq\Phi_L.
\]
Thus $wr_bw^{-1}\in W_L$, and so $r_b\in w^{-1}W_Lw=F$. Since $F$ was an arbitrary
maximal finite subgroup of~$W$ containing~$r_a$, this shows that $r_b\in\FC(r_a)$.

Note that since $W$ has a graph automorphism that swaps $r_a$ and $r_b$ and
fixes all the other simple reflections, it must also be true that $r_a\in\FC(r_b)$.

Let $\widetilde M$ be the component of $\EOdd(a)$ containing $\Odd(a)$, and suppose,
for a contradiction, that $\widetilde M$ is spherical. Note that
$\widetilde M\ne\Odd(a)$, in view of condition~(5) of Definition~\ref{focus}. So
$\widetilde M\setminus\Odd(a)\ne\emptyset$. But for all
$e\in\widetilde M\setminus\Odd(a)$ and $c\in\Odd(a)$ we have
$m_{ce}\ne\infty$, since $\widetilde M$ is spherical, and by conditions
(1)~and~(5) of Definition~\ref{focus} it follows that $m_{ce}=2$ for all
$c\in\Odd(a)$. This contradicts the fact that $\widetilde M$ is connected.

Suppose that $a'\in\Odd(a)$ is such that $\FC(r_{a'})=W_J$ for
some~$J\subseteq\Pi$, and let $J_0$ be the component of $J$ containing~$a'$.
Since $J_0\ne\widetilde M$ it follows from Proposition~\ref{components}
that $J\cap\Odd(a)$ has rank at most~2. Now suppose, for a contradiction, that
$a'\notin\{a,b\}$. Since there exists an element $w\in W_{D[a,b..a']}$ such that
$wa=a'$ and $wb=X(a')$, it follows that
$r_{wb}=wr_bw^{-1}\in\FC(wr_aw^{-1})=\FC(r_{a'})$. Thus $X(a')\in\Phi_J$,
and so $D[a,b..a']\subseteq J$, contradicting the fact that the rank
of $J\cap\Odd(a)$ is at most~2. So we deduce that $a'=b$ or $a'=a$.
Moreover, in either case we know that $\{a,b\}\subseteq J\cap\Odd(a)$,
and since $J\cap\Odd(a)$ has rank at most~2 it follows that
$J\cap\Odd(a)=\{a,b\}$.

By Lemma~\ref{lem1a} we know that all spherical components of $\EOdd(a)$
are components of~$J$, and by
Proposition~\ref{components} all other components of $J$ that intersect $\Odd(a)$
trivially correspond to $C_3$-neighbours of~$\Odd(a)$. But clearly the
conditions of Definition~\ref{halffocus}
imply that $\Odd(a)$ has no $C_3$-neighbours. So we we conclude that $J$
is the union of $\{a,b\}$ and the spherical components of~$\EOdd(a)$,
as required.
\end{proof}

\subsection*{Proof of Theorem~\ref{main}}.

Let $M$ be a connected component of $\oddgraph(\Pi)$, and write
$\widetilde M$ for the component of $\Even(M)$ containing~$M$.

Suppose first that $\widetilde M$ is spherical, so that the conditions of
Case~A of Theorem~\ref{main} are satisfied, and let $a\in M$ be arbitrary.
Observe that all $C_3$-neighbours of $M$ are contained in~$\widetilde M$. Choose
$a'\in M$ such that
$\FC(r_{a'})$ is visible, and let $\FC(r_{a'})=W_J$. By Lemma~\ref{lem1a}
we know that $\widetilde M$ is contained in~$J$, and hence $a\in J$. So by
Proposition~\ref{finclos} it follows that $\FC(r_a)=W_J$ also.
By Proposition~\ref{components} the only possible components of $J$
apart from $J_0$ are the other spherical components of $\Even(M)$, and
by Lemma~\ref{lem1a} all of these are indeed components of~$J$. So $J$
consists of the spherical components of $\Even(M)$, as required.

Now suppose that $\widetilde M$ is not spherical. If there exists
a $b\in\Pi\setminus M$ such that $(a,b)$ is a focus of~$M$ then
it follows from Proposition~\ref{alt2+} that $\FC(r_a)=W_J$, where
$J$ is composed of $\{a,b\}$ and the spherical components of~$\Even(M)$,
and $\FC(r_{a'})$ is not visible for any $a'\in\Odd(a)\setminus\{a\}$.
Similarly, if there exists a $b\in M$ such that $\{a,b\}$ is a half-focus
of~$M$, then it follows from Proposition~\ref{alt3+} that
$\FC(r_a)=\FC(r_b)=W_J$, where $J$ is composed
of $\{a,b\}$ and the spherical components of~$\Even(M)$, and
$\FC(r_{a'})$ is not visible for any $a'\in\Odd(a)\setminus\{a,b\}$.

Finally, suppose that $\widetilde M$ is not spherical and $M$ does not
have a focus or a half focus. Suppose that $a\in M$ is such that
$\FC(r_a)=W_J$ for some $J\subseteq\Pi$, and let $K$ be the component
of~$J$ containing~$a$.

Suppose first that alternative (b) of Proposition~\ref{components}
holds, so that $K=\{a,b\}$ is of type~$C_2$, and $J\cap M=\{a\}$.
Since $M$ does not have any focus in~$\Pi$, it follows from
Proposition~\ref{alt2} that $M\cup\{b\}$ is a spherical component
of~$\Even(M)$. But the component of $\Even(M)$ containing $M$
is $\widetilde M$, which, by our assumptions, is not spherical. So this
case does not arise.

Alternative~(c) of Proposition~\ref{components} is similarly impossible,
by Proposition~\ref{alt3}, and alternative~(e) is also incompatible with
the assumption that~$\widetilde M$ is not spherical. So we conclude
that alternative~(a) holds: $K=\{a\}=J\cap\Odd(a)$. Note that
all spherical components of $\Even(M)$ are components of~$J$, and by
Proposition~\ref{components} the only other possible components are
the sets $\{b\}$ such that $b$ is a $C_3$-neighbour of~$M$.

Suppose that $b$ is a $C_3$-neighbour of $M$ that is adjacent to~$a$.
Let $\tilde a$ be the unique neighbour of~$a$ in~$M$. By Lemma~\ref{morec3}
we know that $r_b\in\FC(r_{\tilde a})$, and so it follows that
$\FC(r_{a})=r_{\tilde a}r_{a}\FC(r_{\tilde a})r_{a}r_{\tilde a}$
contains the reflection along the
root~$r_{\tilde a}r_ab=b+\sqrt 2a+\sqrt 2\tilde a$. Since
$\FC(r_{a})=W_J$, it follows that both $b$ and $\tilde a$ are in~$J$.
But this contradicts the fact that the component of $J$ containing~$a$
is just~$\{a\}$.

This reasoning has shown that $a\in M$ is adjacent in $\Pi$ to
a $C_3$-neighbour of~$M$ then $\FC(r_a)$ is not visible. On the other
hand, we know that there is at least one $a\in M$ such that $\FC(r_a)$
is visible. So we may choose an $a\in M$ such that $\FC(r_a)=W_J$
for some $J\subseteq\Pi$. Since $a$ is not adjacent to
any~$C_3$-neighbour of~$M$ it follows by Lemma~\ref{morec3} that
all $C_3$-neighbours of $M$ are in~$J$. So we conclude that
$J=J'\cup\{a\}$, where $J'$ is the union of the spherical components
of $\Even(M)$ and the $C_3$-neighbours of~$M$.

It remains to prove that is $a'$ is any other element of $M$ that is
not adjacent to any $C_3$-neighbour of~$M$ then
$\FC(r_{a'})=W_{J'\cup\{a'\}}$. Given such an~$a'$, since $a'$ lies in
$M=\Odd(a)$, we may choose $w\in W$ such that $wa=a'$. By
Proposition~\ref{normalizers} we see that
$w\in W_{\widetilde M}$, and so $w$ fixes all other components
of~$\Even(M)$. And $w$ fixes all $C_3$-neighbours
of~$M$, by Lemma~\ref{c3lem}. So $w$ fixes~$J'$, and it follows that
\[
\FC(r_{a'})=w\FC(r_a)w^{-1}=wW_{J'\cup\{a\}}w^{-1}=W_{wJ'\cup\{wa\}}=
W_{J'\cup\{a'\}},
\]
as required. This completes the proof of Theorem \ref{main}.

\subsection*{Proof of Theorem \ref{thm2sph}}

Let $a \in \Pi$ and $M = \Odd(a)$. As $\Pi$ is 2-spherical it follows
that $\Pi = E(M)$, and, as $\Pi$ is non-spherical, it follows that
Case A of Theorem \ref{main} does not hold for $M$.
As there are no $\infty$-labels in the Coxeter graph of $\Pi$, Cases
C and D do not hold either. Hence we are in Case B. Since
there are no $\infty$-labels in the Coxeter graph of $\Pi$, there
are no $C_3$ neighbors of $M$. As $E(M) = \Pi$ is irreducible, there
are no spherical components of $E(M)$. It follows now from 
Theorem \ref{main} that there is an $a' \in \Odd(a)$ such that
$\FC(r_{a'}) = \langle r_{a'} \rangle$. As $r_a$ and $r_{a'}$
are $W$-conjugate we have $\FC(r_a) = \langle r_a \rangle$ as well,
and this completes the proof of Part a) of Theorem~\ref{thm2sph}.

Let $S \subseteq W$ be such that $(W,S)$ is a Coxeter system. It follows
from Part a) and Corollary \ref{FCcrit1} that $r_a \in S^W$ for each
$a \in \Pi$, and hence $\{ r_a \mid a \in \Pi \} \subseteq S^W$.
As $\Pi$ is assumed to be non-spherical, irreducible and 2-spherical, it
follows now from the main result of \cite{CM} that there is
an element $w \in W$ such that $\{ r_a \mid a \in \Pi \} = S^w$.
This completes the proof of Part b) of Theorem \ref{thm2sph}.
As Part c) is an immediate consequence of Part b)  
we are done.

\bibliographystyle{amssort}

\end{document}